\documentclass[titlepage,draft,12pt]{article} 
\usepackage{amsfonts,latexsym,pstricks,pst-plot,graphicx} 
\usepackage[a4paper]{geometry}
\newcommand{\ep}{\varepsilon}

\newcommand{\rep}{r_{\ep}}
\newcommand{\roep}{\rho_{\ep}}
\newcommand{\uep}{u_{\ep}}

\newcommand{\tetep}{\theta_{\ep}}

\newcommand{\m}[1]{\au{#1}^{2\gamma}}
\newcommand{\au}[1]{|A^{1/2}#1|}
\newcommand{\auq}[1]{|A^{1/2}#1|^{2}}
\newcommand{\dat}{D(A^{3/2})}
\newcommand{\da}{D(A)}
\newcommand{\dau}{D(A^{1/2})}
\newcommand{\prf}{{\sc Proof.}$\;$}
\newcommand{\qed}{{\penalty 10000\mbox{$\quad\Box$}}}
\newtheorem{thm}{Theorem}[section]
\newtheorem{thmbibl}{Theorem}
\newtheorem{rmk}[thm]{Remark}
\newtheorem{prop}[thm]{Proposition}

\newtheorem{lemma}[thm]{Lemma}

\title{Hyperbolic-parabolic singular perturbation for mildly 
degenerate Kirchhoff  equations with weak dissipation}
\author{Marina Ghisi\medskip\\ {\normalsize
Universit\`a degli Studi di Pisa} \\{\normalsize Dipartimento di
Matematica ``Leonida Tonelli''}\\
{\normalsize 
PISA (Italy)}\\  
{\normalsize e-mail: \texttt{ghisi@dm.unipi.it}}}
\date{}

\begin{document}
\maketitle
\begin{abstract}
	We consider Kirchhoff equations with a small parameter $\ep$  such as
$$
	\ep\uep''(t)+(1+t)^{-p}\uep'(t)+
	\m{\uep(t)}A\uep(t)=0.
$$
	We prove the existence of global solutions when $\ep$ is 
	small with respect to the size of initial data, for all $0\leq p \leq 
	1$ and $\gamma \geq 1$. Then we provide global-in-time error 
	estimates  on $\uep - u$ where $u$ is the solution of the 
	parabolic problem obtained setting formally $\ep = 0$ in the 
	previous equation.  
	
\vspace{1cm}

\noindent{\bf Mathematics Subject Classification 2000 (MSC2000):}
35B25, 35L70, 35B40.

\vspace{1cm} 

\end{abstract}
 
\section{Introduction}

Let $H$ be a real Hilbert space.  For every $x$ and $y$ in
$H$, $|x|$ denotes the norm of $x$, and $\langle x,y\rangle$ denotes
the scalar product of $x$ and $y$.  Let $A$ be a self-adjoint linear
operator on $H$ with dense domain $D(A)$.  We assume that $A$ is
nonnegative, namely $\langle Ax,x\rangle\geq 0$ for every $x\in D(A)$,
so that for every $\alpha\geq 0$ the power $A^{\alpha}x$ is defined
provided that $x$ lies in a suitable domain $D(A^{\alpha})$.

 For every $\ep>0$ we consider the 
Cauchy problem
\begin{equation}
	\ep\uep''(t)+(1+t)^{-p}\uep'(t)+
	\m{\uep(t)}A\uep(t)=0,
	\label{pbm:h-eq}
\end{equation}
\begin{equation}
	\uep(0)=u_0,\hspace{3em}\uep'(0)=u_1.
	\label{pbm:h-data}
\end{equation}
Equation (\ref{pbm:h-eq}) is the prototype of all degenerate Kirchhoff
equations with weak dissipation, such as
\begin{equation}
	\ep\uep''(t)+(1+t)^{-p}\uep'(t)+
	m(\auq{\uep(t)})A\uep(t)=0 \hspace{2em}
	\forall t\geq 0,
	\label{pbm:h-eq-gen}
\end{equation}
where  $m:[0,+\infty[\to [0,+\infty[$
is a given function which is always assumed to be of class $C^{1}$.  It is
well known that (\ref{pbm:h-eq-gen}) is the abstract setting of a
quasilinear nonlocal partial differential equation of hyperbolic type
which was proposed as a model for small vibrations of strings and
membranes.

Let us start by recalling some general terminology. Equation (\ref{pbm:h-eq-gen}) 
is called \emph{nondegenerate} (or strictly
hyperbolic) when
$$\mu:=\inf_{\sigma\geq 0}m(\sigma)>0,$$
and \emph{mildly degenerate} when $\mu=0$ but $m(\auq{u_{0}})\neq 0$. In the
special case of equation (\ref{pbm:h-eq}) this assumption reduces to 
\begin{equation}
	A^{1/2}u_{0}\neq 0.
	\label{hp:mdg}
\end{equation}
Concerning the dissipation term, we have \emph{constant
dissipation} when $p=0$  and \emph{weak
dissipation} when $p> 0$.  Finally, the operator
$A$ is called \emph{coercive} when
$$	\nu:=\inf\left\{\frac{\langle Ax,x\rangle}{|x|^{2}}:x\in D(A),\
	x\neq 0\right\}>0,$$
and \emph{noncoercive} when $\nu=0$.

In the following we recall briefly what is a singular perturbation 
problem and the ``state of the art''. For a more complete discussion 
on this argument
we refer to the survey \cite{trieste} and to the references  
contained therein.
Moreover we concentrate mostly on equation (\ref{pbm:h-eq}) recalling 
only a few facts on the general equation (\ref{pbm:h-eq-gen}).

The singular perturbation problem in its generality consists in
proving the convergence of solutions of (\ref{pbm:h-eq}),
(\ref{pbm:h-data}) to solutions of the first order problem
\begin{equation}
	u'(t)+ (1+t)^{p}\m{u(t)}Au(t)=0, \hspace{2em}
	u(0)=u_{0},
	\label{pbm:par}
\end{equation}
obtained setting formally $\ep=0$ in (\ref{pbm:h-eq}), and
omitting the second initial condition in~(\ref{pbm:h-data}).
In the concrete case, equation in (\ref{pbm:par}) is a partial 
differential equation of parabolic type. 
With a little abuse of notation in the following we refer to 
hyperbolic and parabolic problems (or behavior) also in the abstract 
setting of equations (\ref{pbm:h-eq}) and (\ref{pbm:par}).

Following the approach introduced by J.\ L.\
Lions~\cite{lions} in the linear case, one defines the corrector
$\tetep(t)$ as the solution of the second order \emph{linear} problem
\begin{equation}
    \ep\tetep''(t)+(1+t)^{-p}\tetep'(t)=0 \hspace{2em}
	\forall t\geq 0,
    \label{theta-eq}
\end{equation}
\begin{equation}
	\tetep(0)=0,\hspace{2em}\tetep'(0)=u_1+
	\m{u_{0}}Au_{0}=:w_{0}.
	\label{theta-data}
\end{equation}
It is easy to see that $\tetep'(0)=\uep'(0)-u'(0)$, hence this
corrector keeps into account the boundary layer due to the loss of one
initial condition.  Finally one defines $\rep(t)$ and $\roep(t)$ in
such a way that
\begin{equation}
	\uep(t)=u(t)+\tetep(t)+\rep(t)=u(t)+\roep(t)\quad\quad\forall
	t\geq 0.
	\label{defn:rep}
\end{equation}
With these notations, the singular perturbation problem consists in
proving that $\rep(t)\to 0$ or $\roep(t)\to 0$ in some sense as
$\ep\to 0^{+}$.  In particular 
	\emph{time-independent} estimates on $\roep(t)$ or $\rep(t)$ as
	$\ep\to 0^{+}$ are called \emph{global error estimates}. 

In this paper  we
restrict ourself to the so called parabolic regime, namely to the 
case where $0\leq p \leq 1$.  The reason is that equations
(\ref{pbm:h-eq}) and (\ref{pbm:h-eq-gen}) have a different behavior
when $p\leq 1$ or $p>
1$.  This is true also in the linear nondegenerate case.
Let us indeed consider equation
\begin{equation}
	au''(t)+\frac{b}{(1+t)^{p}}u'(t)+cAu(t)=0,
	\label{eq:wirth}
\end{equation}
where $a$, $b$, $c$ are positive parameters, and $p\geq 0$.  This
equation was investigated by T.\ Yamazaki~\cite{yamazaki-lin} and J.\
Wirth~\cite{wirth}.  They proved that (\ref{eq:wirth}) has both
parabolic and hyperbolic features, and which nature prevails depends
on $p$.  When $p<1$ the equation has parabolic behavior, in the sense
that all its solutions decay to 0 as $t\to +\infty$ as solutions of
the parabolic equation with $a=0$.  When $p>1$ the same equation has
hyperbolic behavior, meaning that every solution is asymptotic to a
suitable solution of the non-dissipative equation with $b=0$ (and in
particular all non-zero solutions do not decay to zero).  In the
critical case $p=1$ the nature of the problem depends on $b/a$, with
the parabolic behavior prevailing as soon as the ratio is large
enough. In \cite{{gg:w-ndg}} and \cite{jde2} it was proved that also in the case 
of Kirchhoff equation we have always hyperbolic behavior when $p>1$, meaning
that non-zero global solutions (provided that they exist) cannot decay
to 0.  On the other hand, solutions of the limit parabolic problem decay to zero also for $p>1$,
faster and faster as $p$ grows.

The study of the singular perturbation  problem has generated a considerable literature in 
particular regarding the preliminary problem of the existence of 
global solutions for (\ref{pbm:h-eq}) (or (\ref{pbm:h-eq-gen})). 
Despite of this, existence of global solutions without smallness assumptions on
	$\ep$ is a widely open question. 
	The existence of global solutions for (\ref{pbm:h-eq}), 
	(\ref{pbm:h-data}) in the case of a constant dissipation 
	($p=0$) and $\gamma \geq 1$ when $\ep$ is small and (\ref{hp:mdg}) holds true,
	was established by  K.\ Nishihara and
Y.\ Yamada~\cite{ny} (see also  E.\ De
Brito~\cite{debrito1} and Y.\ Yamada~\cite{yamada} for the 
nondegenerate case, \cite{gg:k-dissipative}  for 
the general case and  \cite{ghisi2} for the case $\gamma < 1$).
Moreover optimal and $\ep$-independent decay
estimates  were  obtained in 
\cite{gg:k-decay} (and by  T.\ Mizumachi
(\cite{mizu-ade,mizu-nc}) and K.\ Ono
(\cite{ono-kyushu,ono-aa}) when $\gamma = 1$). 

 When $0\leq p \leq 1$ the existence of  global solutions, always for 
 $\ep$ small, in 
the nondegenerate case was proved 
 in recent years
 by  M.\ Nakao and  J.\
Bae~\cite{nakao}, by T.\
Yamazaki~\cite{yamazaki-wd,yamazaki-cwd}, and in~\cite{gg:w-ndg}. 

The first result for (\ref{pbm:h-eq}) when $p> 0$ 
 was obtained by K.\
Ono~\cite{ono-wd}.  In the special case $\gamma=1$ he proved that a
global solution exists provided that $\ep$ is small and
$p\in[0,1/3]$. Then for ten years there were no significant progresses.
The reason of the slow progress in this field is
hardly surprising. In the weakly dissipative case existence and decay
estimates have to be proved in the same time. The better are the decay
estimates, the stronger is the existence result. This is due to the 
competition between the smallness of the dissipation term and the one 
of the nonlinear term. Both of them decay to zero at  infinity and 
it seems fundamental to understand which of them prevails.
Ten years ago decay estimates for degenerate equations were far from
being optimal, but for the special case $\gamma=1$.
In~\cite{gg:k-decay} a new method for obtaining optimal decay
estimates was introduced and this allowed a substantial progress. In
particular in \cite{jde2} the following result has been proved.
\begin{thmbibl} [\cite{jde2}] \label{A} \em{  Let $(u_{0},u_{1})\in D(A)\times D(A^{1/2})$. If the operator $A$ is coercive and $0\leq 
    p \leq 1$, $\gamma > 0$,
for $\ep$ small the mildly degenerate  problem (\ref{pbm:h-eq}),
        (\ref{pbm:h-data}) has a unique global solution  such that
        $$\frac{C_{1}}{(1+t)^{(p+1)/\gamma}}\leq
                \auq{\uep(t)}\leq
                \frac{C_{2}}{(1+t)^{(p+1)/\gamma}}
                \quad\quad\forall t\geq 0,$$
        $$\frac{C_{1}}{(1+t)^{(p+1)/\gamma}}\leq
                |A\uep(t)|^{2}\leq
                \frac{C_{2}}{(1+t)^{(p+1)/\gamma}}
                \quad\quad\forall t\geq 0,$$
        $$|\uep'(t)|^{2}\leq
                \frac{C_{2}}{(1+t)^{2+(p+1)/\gamma}}
                \quad\quad\forall t\geq 0.$$
		If the operator $A$ is only non negative, $\gamma \geq 1$  and
        \begin{equation}
                0\leq p\leq \frac{\gamma^{2}+1}{\gamma^{2}+2\gamma-1},
                \label{hp:pnc}
        \end{equation}
     for $\ep$ small the mildly degenerate  problem (\ref{pbm:h-eq}),
        (\ref{pbm:h-data}) has a unique global solution  such that
        $$\frac{C_{1}}{(1+t)^{(p+1)/\gamma}}\leq
                \auq{\uep(t)}\leq
                \frac{C_{2}}{(1+t)^{(p+1)/(\gamma+1)}}
                \quad\quad\forall t\geq 0,$$
        $$|A\uep(t)|^{2}\leq
                \frac{C_{2}}{(1+t)^{(p+1)/\gamma}}
                \quad\quad\forall t\geq 0,$$
        $$|\uep'(t)|^{2}\leq
                \frac{C_{2}}{(1+t)^{[2\gamma^{2}+(1-p)\gamma+p+1]/
                (\gamma^{2}+\gamma)}}
                \quad\quad\forall t\geq 0.$$}
		\end{thmbibl}
In the case of noncoercive operators this result is not optimal 
because of (\ref{hp:pnc}). This 
gap is due to the fact that in this second case the estimate on 
$\auq{\uep}$ is worse. This problem is in some sense unavoidable. 
Indeed,
also in the case of linear equations, small 
eigenvalues can make worse the decay of solutions.  Despite of this,
the first result in 
this paper fills up the gap. The key technical point is that in the 
noncoercive case a better decay of $|A\uep|^{2}$ compensates a worse 
decay of $\auq{\uep}$  and the two decay rates are strictly related
(see Proposition \ref{priori}).
This unexpected decay property requires some new and subtle 
estimates. Such estimates improve the decay rates also
when (\ref{hp:pnc})  is satisfied.

Once we know that a global solution of (\ref{pbm:h-eq}) exists, we can 
focus on the singular perturbation problem. This question was solved in the nondegenerate 
case. In that case decay - error estimates were proved,
that consist in estimating in the same time the behavior of $\uep(t)- 
u(t)$ as
	$t\to+\infty$ and as $\ep\to 0^{+}$  (see H.\ Hashimoto and T.\
Yamazaki~\cite{yamazaki}, T.\
Yamazaki~\cite{yamazaki-wd,yamazaki-cwd} and  ~\cite{gg:w-ndg}). 
On the contrary the singular perturbation problem is still quite open in the
degenerate case. With respect to global in time error estimates  we 
indeed know  only the following partial result in
the constant dissipative case (see~\cite{gg:k-PS} where however more 
general nonlinearities are considered).

\begin{thmbibl}[Constant dissipation, \cite{gg:k-PS}] \em{
        \label{thm:error}
        If we assume that $p=0$, $\gamma \geq 1$, 
                $(u_{0},u_{1})\in\dat\times\dau$, then there exists a
                constant $C$ such that for every $\ep$ small we have
                that
                $$|\roep(t)|^{2}+\ep|A^{1/2}\roep(t)|^{2}\leq C\ep^{2}
                \quad\quad
                \forall t\geq 0,$$
                $$\int_{0}^{+\infty}
                |\rep'(t)|^{2}\,dt\leq C\ep.$$}
\end{thmbibl}
This result is far from being optimal.  First of all
only the convergence rate of $|\roep(t)|^{2}$ is optimal, while with 
these regularity assumptions on the initial data one can expect an 
optimal convergence rate also for $\auq{\roep}$ (see 
\cite{gg:l-cattaneo}). Moreover it is limited to equations
with constant dissipation. The second result of this paper fills up 
completely this gap with respect to error estimates in the case of coercive
operators or when $p$ 
verifies (\ref{hp:pnc}) (hence always in the case of a constant 
dissipation) and provides global, but not optimal,  error 
estimates in the remaining cases.  Also to prove this second result 
are fundamental decay estimates as accurate as possible. Conversely 
it is still open the problem of decay-error estimates in all 
degenerate cases.

This paper is organized as follows.  In Section \ref{stat} we state 
precisely our results. Section \ref{proofs} is devoted to the proofs, 
and it is divided into several parts. 
In particular to begin with in Section \ref{sec:ode} 
we state and prove some general lemmata, then in Section \ref{parab}
we consider the parabolic problem (\ref{pbm:par})  and finally  in 
Sections  \ref{ex} and \ref{error} we prove 
the results.

\setcounter{equation}{0}
\section{Statements} \label{stat}
The first result we state concerns the existence of global solutions 
for (\ref{pbm:h-eq}) and their decay properties.
\begin{thm}[Global solutions and decay] \label{th:ex}
    Let us assume that $0 \leq p \leq 1$, $\gamma \geq 1$ and $A$ be 
    a nonnegative operator. Let us assume that $(u_{0}, u_{1}) \in 
    D(A)\times D(A^{1/2})$ satisfy (\ref{hp:mdg}).
    
	Then there exists $\ep_{0}>0$ such that for every
	$\ep\in(0,\ep_{0})$ problem (\ref{pbm:h-eq}), (\ref{pbm:h-data})
	has a unique global solution 
$$
		\uep\in C^{2}([0,+\infty[;H)\cap C^{1}([0,+\infty[;\dau) \cap
		C^{0}([0,+\infty[;\da).
		$$
		Moreover there exist positive constants $C_{1}$ and $C_{2}$ such that
\begin{equation}
    |\uep(t)|^{2} \leq C_{1}  \quad\quad\forall t\geq 0;
    \label{D0}
\end{equation}
	\begin{equation}
		\frac{C_{1}}{(1+t)^{(p+1)/\gamma}}\leq
		\auq{\uep(t)}\leq
		\frac{C_{2}}{(1+t)^{(p+1)/(\gamma+1)}}
		\quad\quad\forall t\geq 0;
		\label{D1}
	\end{equation}
	\begin{equation}
	\au{\uep(t)}^{2(\gamma-1)}	|A\uep(t)|^{2}\leq
		\frac{C_{2}}{(1+t)^{(p+1)}}
		\quad\quad\forall t\geq 0;
		\label{D2}
	\end{equation}
	\begin{equation}
		|\uep'(t)|^{2}\leq 
		\frac{C_{2}	\au{\uep(t)}^{2(\gamma+1)}}{(1+t)^{1-p}}
		\quad\quad\forall t\geq 0;
		\label{D3}
	\end{equation}
	\begin{equation}
	\int_{0}^{+\infty}|\uep'(t)|^{2} (1+t)\, dt \leq 
		C_{2}.
		\label{D4}
	\end{equation}
 \end{thm}
    \begin{rmk}
    \emph{Inequality (\ref{D2}) as far we know is new and it is the core of
    the existence theorem.
    This estimate says that in the case when the operator $A$ is noncoercive 
    it is of course possible that $\au{\uep}$ decays  slow, but in this case 
    $|A\uep|$ decays stronger than in the coercive case.}
    \end{rmk}
    \begin{rmk}
    \emph{When the operator $A$ is noncoercive, $\gamma = 1$ is the only case in 
    which this
    result is contained in  Theorem \ref{A} 
    since in the other cases $(\gamma^{2}+1)/(\gamma^{2} + 2\gamma - 
    1) < 1$ and moreover there is always an improvement on the 
    decay rates (see (\ref{D2}) and (\ref{D3})).}
    \end{rmk}
    The next result regards error estimates. It is divided into two 
    parts. The first one concerns all non negative  operators  and all 
    $0\leq p \leq 1$ and the exponent of $\ep$ in the estimates is 
    not optimal. The second one gives optimal estimates but with some 
    restrictions on the operator or on the admissible values of $p$.
    \begin{thm}[Global-in-time error estimates]
        \label{thm:dg-error}
         Let us assume that $0 \leq p \leq 1$, $\gamma \geq 1$ and $A$ be 
    a nonnegative operator. Let $\uep(t)$
        be the solution of equation (\ref{pbm:h-eq}) with 
	initial data $(u_{0},u_{1})\in D(A^{3/2})\times\dau$ 
	satisfying
        (\ref{hp:mdg}).  Let $u(t)$ be the solution of the corresponding
       first order problem (\ref{pbm:par}), and let $\rep(t)$ and $\roep(t)$ be defined by
        (\ref{defn:rep}).

        Then we have the following conclusions.
        \begin{enumerate}
                \renewcommand{\labelenumi}{(\arabic{enumi})}
                \item There exists a
                constant $C_{3}$ such that for every $\ep$ small 
		enough we have
                that
                \begin{equation}
                    |\roep(t)|^{2}+|A^{1/2}\roep(t)|^{2} + \int_{0}^{t}
                \frac{|\rep'(s)|^{2}}{\au{u(s)}^{2\gamma}}\frac{1}{(1+s)^{p}}ds 
		\leq C_{3}\ep
                \quad\quad
                \forall t\geq 0.
                    \label{H1}
                \end{equation}
                     \item  If in addition we assume that 
               	\begin{equation}
		    \mbox{$A$ is coercive or } \; 0 \leq p \leq 
		\frac{\gamma^{2}+1}{\gamma^{2}+2\gamma - 1}
		    \label{H1cip}
		\end{equation}
	 then there exists a
                constant $C_{4}$ such that for every $\ep$ small 
		enough we have
                that
                \begin{equation}
                    |\roep(t)|^{2}+|A^{1/2}\roep(t)|^{2} + \int_{0}^{t}
                \frac{|\rep'(s)|^{2}}{\au{u(s)}^{2\gamma}}
		\frac{1}{(1+s)^{p}}ds \leq C_{4}\ep^{2}
                \quad\quad
                \forall t\geq 0.
                    \label{H1c}
                \end{equation}
        \end{enumerate}
\end{thm}
  \begin{rmk} \emph{When the initial data are more regular it is of 
      course possible to 
      achieve an estimate on $A\roep$ like the ones in 
      (\ref{H1}) and (\ref{H1c}). Moreover in this case one can get 
      also estimates on  $\rep'$ exactly as in \cite{gg:k-PS}, \cite{gg:w-ndg} (see also 
      \cite{trieste})  We do not give here the precise 
      statements and proofs since they only lengthen the paper  
      without introducing new ideas.}
 \end{rmk}
 \begin{rmk} \emph{In the integrals in (\ref{H1}), (\ref{H1c}) it appears the 
     coefficient $\au{u}^{-2\gamma}$. When $A$ is a coercive operator 
     we can replace this term with 
      $\au{\uep}^{-2\gamma}$  or $(1+t)^{p+1}$, indeed they 
     all have the same behavior. On the contrary when $A$ is noncoercive 
     the use of $\au{u}^{-2\gamma}$ seems compulsory.}
 \end{rmk}
\setcounter{equation}{0}
\section{Proofs}\label{proofs}
Proofs are organized as follows. First of all in Section 
\ref{sec:ode} we state and prove some general lemmata that do not 
concern directly the Kirchhoff equation. In Section \ref{parab} then 
we recollect all the properties of the solutions of (\ref{pbm:par}) we 
need. Finally in Section \ref{ex} we prove  Theorem \ref{th:ex} and 
in Section \ref{error} we prove Theorem \ref{thm:dg-error}.
\subsection{Basic Lemmata}\label{sec:ode}

Numerous variants of the following comparison result have already been
used in
\cite{ghisi2,gg:k-dissipative,gg:k-decay,gg:w-ndg, 
jde2} and we refer to these ones for the proof.

\begin{lemma}\label{lemma:ode}
	Let $T>0$,  and let $f:[0,T[\to[0,+\infty[$ be a
	function of class $C^{1}$. Let $\phi:[0,T[\to[0,+\infty[$ be 
	a  continuous function.  Then the following implications hold 
	true.  
	
\begin{enumerate}
		\renewcommand{\labelenumi}{(\arabic{enumi})} 
		\item\label{ODEA1}	Let us assume that there exists 
	a constant $a\geq 0$  such that
	$$ f'(t)\leq 
		-\phi(t)\sqrt{f(t)}(\sqrt{f(t)}-a)
		\hspace{2em}
		\forall t\in [0,T[;
	$$
	then we have that
	$$
		f(t)\leq \max\{f(0), a^{2}\} 
		\hspace{2em}
		\forall t\in [0,T[.
	$$
	\item \label{ODEA2}
Let us assume that  there exists 
	a constant $a\geq 0$  such that
	$$ f'(t)\leq 
		-\phi(t)f(t)(f(t)-a)
		\hspace{2em}
		\forall t\in [0,T[;
	$$
	then we have that
	$$
		f(t)\leq \max\{f(0), a\} 
		\hspace{2em}
		\forall t\in [0,T[.
	$$
		\end{enumerate}
\end{lemma}

A proof of the next comparison result  is contained in \cite{jde2} 
(Lemma 3.2).
\begin{lemma}\label{lemma:ode-vs}
 Let $w:[0,+\infty[\to[0,+\infty[$ be a function of class $C^{1}$ with
 $w(0)>0$. Let $a>0$ be a positive constant.

	Then the following implications hold true.
	\begin{enumerate}
		\renewcommand{\labelenumi}{(\arabic{enumi})} \item
		\label{ODEB2} If
		$w$ satisfies the differential inequality $$
			w'(t)\leq
			-a(1+t)^{p}\left[w(t)\right]^{1+\gamma}
			\quad\quad\forall t\in[0,+\infty[, $$
		then for some constant $\gamma_{1}$ we have the
		following estimate
	$$
			w(t)\leq
			\frac{\gamma_{1}}{(1+t)^{(p+1)/\gamma}}
			\quad\quad\forall t\in[0,+\infty[.
			$$
		\item \label{ODEB1}If $w$ satisfies the differential inequality $$
			w'(t)\geq
			-a(1+t)^{p}\left[w(t)\right]^{1+\gamma}
			\quad\quad\forall t\in[0,+\infty[, $$
		then for some constant $\gamma_{2}$ we have the
		following estimate
		$$
			w(t)\geq
			\frac{\gamma_{2}}{(1+t)^{(p+1)/\gamma}}
			\quad\quad\forall t\in[0,+\infty[.
		$$
	\end{enumerate}
\end{lemma}
Let us now state and prove the third lemma.
\begin{lemma}\label{lemma:ode-GF}
 Let $F,\, G:[0,T[\to[0,+\infty[$ be functions of class $C^{1}$.  Let
 $\varphi:[0,T[\to]0,+\infty[$ be a continuous function and $a>0$,
 $b \geq 0,\, c \geq 0$ be real numbers.  Let us assume that in $[0,T[$ the 
 following inequality holds true:
\begin{equation}
    (F+G)'(t) \leq - \varphi(t) (F(t) + a (G(t))^{2} -bG(t) -
    c(G(t))^{3/2}).
    \label{0-ode}
\end{equation}
Let us set $\sigma_{0}:= (c + \sqrt{ab})/a$.  Then we get
\begin{equation}
    G(t) + F(t) \leq \sigma_{0}^{2}(1 + b + c\sigma_{0}) + F(0) + G(0)
    + 1, \quad\quad \forall t\in[0,T[.
    \label{1-ode}
\end{equation}
\end{lemma}
\prf Let us set
$$S := \sup\{t< T: \:\mbox{ in $[0,t]$ the inequality (\ref{1-ode})
    holds true}\}.
    $$
It is obvious that $S > 0$.  We want to prove that $S = T$.  Let us
assume by contradiction that $S < T$.  Therefore in $[0,S[$ the
inequality (\ref{1-ode}) holds true, moreover
\begin{equation}
     G(S) + F(S) = \sigma_{0}^{2}(1 + b + c\sigma_{0}) + F(0) +
     G(0) + 1
    \label{3-ode}
\end{equation}
and
\begin{equation}
    (G+F)'(S) \geq 0.
    \label{4-ode}
\end{equation}
If $G(S) > \sigma_{0}^{2}$ then 
\begin{equation} a (G(S))^{2}
-bG(S) - c(G(S))^{3/2} = G(S)(aG(S) - b - \sqrt{G(S)} ) > 0.
\label{4b-ode}
\end{equation}
    Indeed let us set $y =\sqrt{G(S)}$, then  
    $$a y^{2} - cy - b > 0 \hspace{1em} \mbox{if}
    \hspace{1em} y > \frac{c + \sqrt{c^{2}+ 4ab}}{2a} =: \sigma_{1}$$
    and  by assumption $\sigma_{1} \leq \sigma_{0} < y$.  Plugging (\ref{4b-ode}) in
    (\ref{0-ode}) we hence arrive at 
    $$ (G+F)'(S) < 0$$
    in contrast with (\ref{4-ode}).
   
   Let us now assume that $G(S) \leq \sigma_{0}^{2}$. Hence from 
   (\ref{3-ode}) we get:
   \begin{eqnarray*}
      & F(S) - bG(S) - c(G(S))^{3/2}  =  & \\
      &\sigma_{0}^{2}(1 + b + c\sigma_{0}) + F(0) +
     G(0) + 1 - (b+1)G(S)  - c(G(S))^{3/2} & \\
        & >   \sigma_{0}^{2}(1 + b + c\sigma_{0}) - (b+1)G(S)  - c(G(S))^{3/2} & \\
        & \geq  \sigma_{0}^{2}(1 + b + c\sigma_{0}) - 
	(b+1)\sigma_{0}^{2}  - c\sigma_{0}^{3} = 0.&
   \end{eqnarray*}
   Hence by (\ref{0-ode}) we obtain once again
 $$ (G+F)'(S) < 0$$
    in contrast with (\ref{4-ode}).
\qed 

 The following lemma is essential in the proof of error estimates.
 \begin{lemma}
     \label{mprop}
     Let us assume that $m:[0,+\infty[\rightarrow[0,+\infty[$ is a 
     nondecreasing function. Then for all $x,\, y \in D(A)$ we get
    $$
         \langle m(\auq{x})Ax - m(\auq{y})Ay, x-y\rangle \geq 
	 \frac{1}{2}\left[m(\auq{x}) + m(\auq{y})\right]\auq{(x-y)}.
         $$
     \end{lemma}
     \prf Let us set
     $$m_{x}:=  m(\auq{x}), \quad\quad m_{y}:= m(\auq{y}).$$
     Thus  an elemental calculation gives:
     \begin{eqnarray*}
        &   \displaystyle \langle m_{x}Ax - m_{y}Ay, x-y\rangle  =   
	    m_{x}\auq{x} +  m_{y}\auq{y}  -(m_{x} + 
	  m_{y})\langle Ax, y\rangle &  \\
          &=\displaystyle  \frac{1}{2}(m_{x} + 
	  m_{y})(\auq{x} + \auq{y}) -  
	   \displaystyle
	 \frac{1}{2}\cdot 2(m_{x} + m_{y})\langle Ax, y\rangle + &   \\
          &  \displaystyle+\frac{1}{2}(m_{x} - m_{y})(\auq{x} - \auq{y})  &  \\
          & =  \displaystyle \frac{1}{2}(m_{x} + 
	  m_{y})\auq{(x-y)} + 
	  \displaystyle \frac{1}{2}(m_{x} - m_{y})(\auq{x} - \auq{y}) &
	   \\
	  &\geq  \displaystyle  \frac{1}{2}(m_{x} + 
	  m_{y})\auq{(x-y)}; &
     \end{eqnarray*}
    where in the last step we exploit that  $m$ is nondecreasing, hence 
    $$(m(\alpha) - m(\beta))(\alpha - \beta) \geq 0 \hspace{1em} 
    \forall \alpha,\, \beta \geq 0.$$
    \qed
    
   The last lemma concerns the integrability properties of the 
   corrector $\tetep$.
   \begin{lemma} \label{thetalemma}
       Let  $0\leq p \leq 1$ and let $\tetep$ be the solution of (\ref{theta-eq}), 
       (\ref{theta-data}). Let $\delta \geq 0$ and let us assume that 
       $\ep<(2+2\delta)^{-1}$. Then  there exists a constant 
       $C_{\delta}$ independent from $\ep$ and from the initial data such 
       that if $w_{0} \in D(A^{j/2})$ therefore we have
       $$\int_{0}^{+\infty}(1+t)^{\delta}|A^{j/2}\tetep'(t)| \, dt \leq 
       C_{\delta}|A^{j/2}w_{0}|\ep.$$
    \end{lemma}
   \prf 
   Let us define
   $$I:= \int_{0}^{+\infty}(1+t)^{\delta}|A^{j/2}\tetep'(t)| \, dt.$$
   If $p=1$ then $\tetep'(t) = w_{0}(1+t)^{-1/\ep}$ hence thesis 
   follows from 
  $$ I = 
  |A^{j/2}w_{0}|\frac{\ep}{1 - (\delta +1)\ep}.$$
  Let us now assume that $p<1$. In such a case we have 
  $$\tetep'(t) = 
  w_{0}\exp\left(-\frac{1}{\ep}\frac{1}{
  1-p}((1+t)^{1-p} - 1)\right).$$ 
  If we set
  $$\phi(t) := \min\{t,t^{1-p}\}$$
  then it is easy to prove that there exists a constant $\beta_{0}> 0$ 
  such that
 $$\frac{1}{1-p}((1+t)^{1-p} - 1)\geq \beta_{0}\phi(t).$$
 
In particular we obtain
\begin{eqnarray*}
     I & \leq & 
     |A^{j/2}w_{0}| 
     \int_{0}^{+\infty}(1+t)^{\delta}e^{ -\frac{1}{\ep}\beta_{0}\phi(t)}dt 
      \\
      & = & |A^{j/2}w_{0}|\left(
     \int_{0}^{1}(1+t)^{\delta}e^{-\frac{t}{\ep}\beta_{0}}dt + 
     \int_{1}^{+\infty}(1+t)^{\delta}e^{-\frac{t^{1-p}}{\ep}\beta_{0}}dt\right). 
 \end{eqnarray*}
    Let us set $\ep s = t$, hence     
 \begin{eqnarray*}
    I & \leq & 
     |A^{j/2}w_{0}| \ep \left(
     \int_{0}^{+\infty}(1+\ep s)^{\delta}e^{-\beta_{0}s}ds + 
     \int_{0}^{+\infty}(1+\ep 
     s)^{\delta}e^{-\frac{s^{1-p}}{\ep^{p}}\beta_{0}}ds\right)
     \\
      & \leq & |A^{j/2}w_{0}|\ep \left(
      \int_{0}^{+\infty}(1+s)^{\delta}e^{-\beta_{0}s}ds + 
     \int_{0}^{+\infty}(1+s)^{\delta}e^{-\beta_{0}s^{1-p}}ds\right)\\
     & = &  |A^{j/2}w_{0}|\ep C_{\delta}.
   \end{eqnarray*}  \qed
   
   \subsection{The First order problem}\label{parab}

Theory of parabolic equations of Kirchhoff type is quite well
established.  These equations appeared for the first time in the
pioneering paper~\cite{bernstein} by S.\ Bernstein and then were 
considered by many
authors (see~\cite{bw,miletta, k-par} and ~\cite{trieste} for the 
details). In fact the following result holds true.

\begin{thmbibl}[Global solutions]\label{thp}
        Let $A$ be a nonnegative operator, let  $0\leq p \leq 1$ and 
	$\gamma \geq 1$.  Let $u_{0}\in\da$.
	
	Then problem (\ref{pbm:par}) has a unique global solution
	$$u\in C^{1}\left([0,+\infty[;H\right)\cap
	C^{0}\left([0,+\infty[;\da\right).$$
	
	If in addition $A^{1/2}u_{0}\neq 0$ then
	the solution is non-stationary, i.e. $\auq{u(t)}\neq 0$ for 
	all $t\geq 0$ and $u\in
	C^{\infty}\left(]0,+\infty[;D(A^{\alpha})\right)$ for every $\alpha\geq
	0$.
\end{thmbibl}

In the  proposition below we collect all the properties of the 
solutions of (\ref{pbm:par})
 we need in  proof of error estimates. Only some of these 
properties require $u_{0} \in \dat$. Nevertheless this is an 
assumption of Theorem \ref{thm:dg-error} hence we do not specify in what 
cases it is in fact necessary or not.

\begin{prop}[Properties of solutions]
Let $u_{0}\in\dat$ and let us assume that all conditions of 
Theorem \ref{thp} are verified.  Let $u$ be the global
solution of (\ref{pbm:par}).
	
	Then the following statements hold true.
	\begin{itemize}
	    \item The solution $u$ verifies the standard estimates 
	    below:
	    \begin{equation}
		\frac{|A^{(k+1)/2}u(t)|^{2}}{|A^{k/2}u(t)|^{2}} \leq
		\frac{|A^{(k+1)/2}u_{0}|^{2}}{|A^{k/2}u_{0}|^{2}}, 
		\quad \quad k=1,\, 2, \quad \quad \forall t\geq 0;
		\label{dv0}
	    \end{equation}
	    \begin{equation}
		\frac{1}{2}|u(t)|^{2}
		+\int_{0}^{t}\au{u(s)}^{2(\gamma+1)}(1+s)^{p} ds =
		\frac{1}{2}|u_{0}|^{2}, \quad \quad \forall t\geq 0.
		\label{dv1}
	    \end{equation}
		    \item The solution $u$ has these decay properties:
	    \begin{equation}
		\frac{\gamma_{3}}{(1+t)^{(p+1)/\gamma}}\leq \auq{u(t)}
		\leq \frac{\gamma_{4}}{(1+t)^{(p+1)/(\gamma+1)}},  
		\quad \quad \forall t\geq 0;
		\label{dv2}
	    \end{equation}
	    \begin{equation}
	\au{u(t)}^{2(\gamma-1)}|Au(t)|^{2} \leq
		\frac{\gamma_{4}}{(1+t)^{p+1}},  \quad \quad \forall t\geq 0.
		\label{dv4}
	    \end{equation}
	If moreover $A$ is a coercive operator then
	    \begin{equation}
		\auq{u(t)} \leq
		\frac{\gamma_{4}}{(1+t)^{(p+1)/\gamma}},  
		\quad \quad \forall t\geq 0.
		\label{dv3}
	    \end{equation}
	     \item The following integrals  are bounded:
	    \begin{equation}
		\int_{0}^{+\infty}|u'(t)|^{2}(1+t)^{p} dt =
		\int_{0}^{+\infty}\au{u(t)}^{4\gamma}|Au(t)|^{2}(1+t)^{3p}
		dt \leq \gamma_{5};
		\label{dv5}
	    \end{equation}
	    \begin{equation}
		\int_{0}^{+\infty}\au{u(t)}^{6\gamma}|A^{3/2}u(t)|^{2}(1+t)^{5p}
		dt \leq \gamma_{5};
		\label{dv6}
	    \end{equation}
	    \begin{equation}
		\int_{0}^{+\infty}\left[\au{u(t)}^{8\gamma}(1+t)^{7p}+
		\au{u(t)}^{6\gamma}(1+t)^{5p}\right]|A^{2}u(t)|^{2} dt \leq 
		\gamma_{5}.
		\label{dv7}
	    \end{equation}
	\end{itemize}
\end{prop}

\prf From now in most of the proofs we omit the dependence of $u$ from
$t$.  Moreover often we use that $0\leq p \leq 1$ but for shortness
sake we do not recall it more. Furthermore we use that in 
$]0,+\infty[$ the solution $u$ is as regular as we want.  

\subparagraph{\textmd{\textsl{Proof of (\ref{dv0})}}} It is enough to
remark that 
$$\left(\frac{|A^{(k+1)/2}u|^{2}}{|A^{k/2}u|^{2}}
\right)' = -2(1+t)^{p}\frac{|A^{1/2}u|^{2\gamma}}{|A^{k/2}u|^{4}}
(|A^{(k+2)/2}u|^{2}|A^{k/2}u|^{2} - |A^{(k+1)/2}u|^{4} ) \leq 0,$$
where in the last inequality we exploit that
\begin{equation}
    |A^{(k+1)/2}u|^{2} =
\langle A^{(k+2)/2}u,A^{k/2}u\rangle \leq |A^{(k+2)/2}u||A^{k/2}u|.
    \label{rappA}
\end{equation}

\subparagraph{\textmd{\textsl{Proof of (\ref{dv1})}}} It suffices 
to integrate in $[0,t]$ the equality: 
$$ \left(\frac{1}{2}|u|^{2} \right)' +\au{u}^{2(\gamma+1)}(1+t)^{p} = 0.$$

\subparagraph{\textmd{\textsl{Proof of (\ref{dv2})}}} Using
(\ref{dv0}) with $k = 1$ we have
$$(\auq{u})'= - 2(1+t)^{p}|Au|^{2}\au{u}^{2\gamma} \geq 
-2(1+t)^{p}
\frac{|A u_{0}|^{2}}{|A^{1/2}u_{0}|^{2}}\au{u}^{2(\gamma+1)} .$$
Therefore  estimate form below follows from Statement (\ref{ODEB1}) in Lemma 
\ref{lemma:ode-vs}.  

Let us now remark that
$$
    \left((1+t)^{p+1}\frac{\au{u}^{2(\gamma +
    1)}}{2(\gamma+1)}\right)'+ (1+t)^{2p+1}\au{u}^{4\gamma}|Au|^{2} =
    \frac{p+1}{2(\gamma+1)}(1+t)^{p}\au{u}^{2(\gamma + 1)}.
    \label{decv0}
$$
Integrating in $[0,t]$ and using (\ref{dv1}) we get:
\begin{equation}
     (1+t)^{p+1}\frac{\au{u(t)}^{2(\gamma + 1)}}{2(\gamma+1)}+
     \int_{0}^{t}(1+s)^{2p+1}\au{u(s)}^{4\gamma}|Au(s)|^{2}ds \leq
     |u_{0}|^{2} + \au{u_{0}}^{2(\gamma+1)}.
    \label{decv}
\end{equation}
From this inequality we gain directly the estimate from above in
(\ref{dv2}).

\subparagraph{\textmd{\textsl{Proof of (\ref{dv4})}}} Let us define
$$G(t) = (1+t)^{p+1}\au{u(t)}^{2(\gamma-1)}|Au(t)|^{2}.$$
We have to prove that $G$ is bounded. Taking the time's derivative of 
$G$ we obtain
$$G' = \frac{-G}{1+t}\left[2\au{u}^{2\gamma}
\frac{|A^{3/2}u|^{2}}{|Au|^{2}}(1+t)^{p+1}+2(\gamma - 1)G
-(p+1)\right].$$
Now let us distinguish two cases.  If $\gamma > 1$ we have:
$$G' \leq\frac{-G}{1+t}\left[2(\gamma - 1)G -(p+1)\right]$$
then thesis follows from Statement (\ref{ODEA2}) in Lemma \ref{lemma:ode}.

Instead if $\gamma = 1$ using (\ref{rappA}) with $k=1$ we get 
$$G'\leq\frac{-G}{1+t}\left[2|Au|^{2}(1+t)^{p+1} -(p+1)\right] =
\frac{-G}{1+t}\left[2G -(p+1)\right], $$
hence we conclude as in the previous case.

\subparagraph{\textmd{\textsl{Proof of (\ref{dv3})}}} Since
$\langle Au, u \rangle \geq \nu |u|^{2}$ then $$(\auq{u})'= -
2(1+t)^{p}|Au|^{2}\au{u}^{2\gamma} \leq -2\nu (1+t)^{p}
\au{u}^{2(\gamma+1)} .$$ Hence it suffices to apply Statement (\ref{ODEB2})
in Lemma \ref{lemma:ode-vs}.

\subparagraph{\textmd{\textsl{Proof of (\ref{dv5})}}} Since $3p\leq 2p+1$ it is a 
consequence of (\ref{decv}).

\subparagraph{\textmd{\textsl{Proof of (\ref{dv6})}}} A simple computation gives:
\begin{eqnarray*}
    & \displaystyle\left(\frac{1}{2}\au{u}^{4\gamma}|Au|^{2}(1+t)^{4p}\right)' +
    \au{u}^{6\gamma}|A^{3/2}u|^{2}(1+t)^{5p} = & \\
   & 2p\au{u}^{4\gamma}|Au|^{2}(1+t)^{3p} 
    -2\gamma\au{u}^{6\gamma-2}|Au|^{4}(1+t)^{5p} \leq
    2\au{u}^{4\gamma}|Au|^{2}(1+t)^{2p+1}.&
\end{eqnarray*}
Hence thesis follows integrating in $[0,t]$ and using (\ref{decv}).

\subparagraph{\textmd{\textsl{Proof of (\ref{dv7})}}} As in the
previous case we have
   \begin{equation}
    \left(\frac{1}{2}\au{u}^{6\gamma}|A^{3/2}u|^{2}(1+t)^{6p}\right)'
    + \au{u}^{8\gamma}|A^{2}u|^{2}(1+t)^{7p} \leq
    3\au{u}^{6\gamma}|A^{3/2}u|^{2}(1+t)^{5p}.
    \label{Iv1}
\end{equation}
    Moreover from (\ref{dv0}) with $k=2$ we get
  \begin{eqnarray}
    &
    \displaystyle \left(\frac{1}{2}\au{u}^{4\gamma}|A^{3/2}u|^{2}(1+t)^{4p}\right)'
    + \au{u}^{6\gamma}|A^{2}u|^{2}(1+t)^{5p} \leq &\nonumber \\
   & 2p\au{u}^{4\gamma}|A^{3/2}u|^{2}(1+t)^{3p} \leq
   2\displaystyle \frac{|A^{3/2}u_{0}|^{2}}{|Au_{0}|^{2}}\au{u}^{4\gamma}|Au|^{2}(1+t)^{3p}.&
    \label{Iv2}
\end{eqnarray}
Summing up (\ref{Iv1}) and (\ref{Iv2}), integrating in $[0,t]$ and using
(\ref{dv6}) and (\ref{dv5}) we end up with (\ref{dv7}).

\subsection{Proof of Theorem \ref{th:ex}}\label{ex}
As in the previous section  in most of the proofs we omit the 
dependence of $\uep$ from 
$t$ and we do not recall more that $0\leq p \leq 1$. We divide the 
proof into three parts. In the first one 
we state and prove the energy estimates we need, then we prove the 
existence of global solutions and finally we give the decay estimates.

\subsubsection{Basic energy estimates}

In this section we prove some estimates that involve the following 
energies:
\begin{eqnarray}
    Q_{\ep}(t) & = & 
    \frac{|\uep'(t)|^{2}}{\au{\uep(t)}^{2(\gamma+1)}}(1+t)^{1-p}; 
    \label{defQ} \\
    D_{\ep}(t) & = & \ep
    \frac{\langle \uep'(t), 
    \uep''(t)\rangle}{\au{\uep(t)}^{2(\gamma+1)}}(1+t)^{2p+1} + \int_{0}^{t} 
    \frac{\auq{\uep'(s)}}{\auq{\uep(s)}}(1+s)^{2p+1}ds; \label{defD}  \\
    R_{\ep}(t) & = & \left[\ep 
    \frac{|\uep''(t)|^{2}}{\au{\uep(t)}^{2(\gamma+1)}} + 
    \frac{\auq{\uep'(t)}}{\auq{\uep(t)}}\right](1+t)^{2(p+1)}; 
    \label{defR}  \\
    H_{\ep}(t) & = & \left[\ep 
   \frac{\auq{\uep'(t)}}{\auq{\uep(t)}} + 
    \au{\uep(t)}^{2(\gamma-1)}|A\uep(t)|^{2}\right](1+t)^{p+1}. 
    \label{defH} 
\end{eqnarray}

Let us moreover set:
$$h_{1}:= 4(|u_{1}|^{2} 
+\au{u_{0}}^{4\gamma}|Au_{0}|^{2})\au{u_{0}}^{-2(\gamma+1)}, \hspace{1em}h_{2}:= (\gamma - 1)(\sqrt{h_{1}} +1) + \sqrt{\gamma - 
1},$$
$$L_{1}:= \left\{
\begin{array}{ll}
    (3 + 2h_{2}(\sqrt{h_{1}} +1))h_{2}^{2}(\gamma - 1)^{-2} + 
    H_{1}(0)+1&  
    \hspace{1em} \mbox{if $\gamma > 1$}  \\
    \\
    36 + 2 \auq{u_{1}}\au{u_{0}}^{-2}+ 2|Au_{0}|^{2} + 2^{-1}|\langle 
    Au_{0}, u_{1}\rangle|\au{u_{0}}^{-2} & \hspace{1em} \mbox{if 
    $\gamma = 1$}.
\end{array}\right.$$
In the following proposition we recollect all the estimates on 
(\ref{defQ}) trough (\ref{defH}) we need.
\begin{prop}[A priori estimates] \label{priori}
  Let us assume that  all the hypotheses of Theorem \ref{th:ex} 
  are verified.  
  Then there exists $\ep_{0}$ with the following property. If
  $\ep\in]0,\ep_{0}]$, $S > 0$ and
  $$\uep\in C^{2}([0,S[;H)\cap C^{1}([0,S[;\dau) \cap
	C^{0}([0,S[;\da)$$
is a solution of (\ref{pbm:h-eq}), (\ref{pbm:h-data}) such that
  \begin{equation}
   \auq{\uep(t)}> 0 \hspace{1em}
    \forall t\in[0,S[,
      \label{fond1}
  \end{equation}
  \begin{equation}
   \frac{|\langle \uep'(t), 
    A\uep(t)\rangle|}{\auq{\uep(t)}} \leq
   \frac{K_{0}}{1+t}, \hspace{1em}
   \au{\uep(t)}^{2(\gamma-1)}|A\uep(t)|^{2} \leq
   \frac{K_{1}}{(1+t)^{p+1}}, \;\; \forall t\in[0,S[,
      \label{fond}
  \end{equation}
  then there exists  a positive constant $L_{3}$ independent from $\ep$ and $S$ such that for every
  $t\in[0,S[$:
	\begin{equation}
	   Q_{\ep}(t)  \leq  \max\{4K_{1}, Q_{1}(0)\}=: L_{2}; \label{SQ}
       \end{equation}
	   \begin{equation} D_{\ep}(t) \leq D_{\ep}(0) + 2 L_{2}(3+2K_{0})
	    (1+t)^{p+1} +\frac{1}{8(K_{0}+1)} \int_{0}^{t}
	    \frac{|\uep''(s)|^{2} (1+s)^{2p+1}}{\au{\uep(s)}^{2(\gamma+1)}}
	   ds; \label{SD}
	   \end{equation}
	   \begin{equation} \left[\ep 
    \frac{|\uep''(t)|^{2}}{\au{\uep(t)}^{2(\gamma+1)}} + 
    \frac{\auq{\uep'(t)}}{\auq{\uep(t)}}\right](1+t)^{p+1}= 
    \frac{R_{\ep}(t)}{(1+t)^{p+1}} \leq L_{3}+ 2R_{\ep}(0); 
	    \label{SR}
	    \end{equation}
	    \begin{equation}  H_{\ep}(t)  \leq  L_{1}.
	    \label{SH}
	\end{equation}
   \end{prop}
\prf
Let us set:
\begin{equation}
    h_{3}:= 
4(4\gamma^{2}K_{0}^{2}K_{1} + L_{2}),\hspace{1em}
h_{4}:=h_{3}+ 8(K_{0}+1)(3+2K_{0})L_{2},
   \label{defsigma1}
\end{equation}
\begin{equation}
    L_{3}:= 2\left[h_{4} + \frac{L_{2}}{2} + 4 
\frac{|u_{1}|}{\au{u_{0}}^{2(\gamma+1)}}(|u_{1}| + 
\au{u_{0}}^{2\gamma}|Au_{0}|)(K_{0}+1)\right].
    \label{defL2}
\end{equation}

Now let us assume that $\ep_{0}$ verifies the following inequalities:
\begin{equation}
      8\ep_{0}(2+(\gamma + 1)K_{0} )\leq 1, \quad\quad 
      16\ep_{0}(K_{0}+1)^{2}\leq 1,
    \label{1-ep}  
    \end{equation}
    \begin{equation}
      2\ep_{0}(K_{0}+1)(1 + (3 + 2(\gamma + 1)K_{0})^{2})\leq 1/8,
    \label{2-ep}
    \end{equation}
    \begin{equation}
  \sqrt{\ep_{0}}\left(\sqrt{L_{3}} +\sqrt{2} 
    \frac{\au{u_{1}}}{\au{u_{0}}}\right) \leq 1.    
    \label{3-ep}
\end{equation}

Let us now compute the time's derivatives of the energies (\ref{defQ}) 
through (\ref{defH}). After some computation we find that:

\begin{equation}
    Q_{\ep}'= Q_{\ep}\left(\frac{(1-p)}{1+t} 
    -2(\gamma+1)\frac{\langle A\uep, \uep'\rangle}{\auq{\uep}}
    -\frac{2}{\ep}\frac{1}{(1+t)^{p}}\right) 
    -\frac{2}{\ep}(1+t)^{1-p}\frac{\langle A\uep, 
    \uep'\rangle}{\auq{\uep}}.
    \label{derQ}
\end{equation}
Let us set:
\begin{eqnarray*}
   \varphi_{1}(t)  &:=  &  \ep 
    \frac{|\uep''(t)|^{2}}{\au{\uep(t)}^{2(\gamma+1)}}(1+t)^{2p+1} +\\
     &  & + \ep\left[  \frac{\langle \uep'(t), 
    \uep''(t)\rangle}{\au{\uep(t)}^{2(\gamma+1)}}(1+t)^{2p}\left(2p+1 
    -2(\gamma+1)\frac{\langle A\uep(t), 
    \uep'(t)\rangle}{\auq{\uep(t)}}(1+t)\right)\right],
\end{eqnarray*}
   $$\varphi_{2}(t) :=- 2\gamma\left(\frac{\langle A\uep(t), 
     \uep'(t)\rangle}{\auq{\uep(t)}}\right)^{2}(1+t)^{2p+1},$$
     $$\varphi_{3}(t):= - \frac{\langle \uep'(t), 
    \uep''(t)\rangle}{\au{\uep(t)}^{2(\gamma+1)}}(1+t)^{p+1},$$
    $$  \varphi_{4}(t):= p \frac{|\uep'(t)|^{2}}{\au{\uep(t)}^{2(\gamma+1)}}
    (1+t)^{p};$$
    thus
\begin{equation}
    D_{\ep}' 
      =  \varphi_{1} +  \varphi_{2} +  \varphi_{3} +  
     \varphi_{4}.
    \label{derD}
\end{equation}
Moreover
\begin{eqnarray}
    R_{\ep}' & = & 2(1+t)^{2(p+1)}\left(-2\gamma \frac{\langle A\uep,
    \uep'\rangle}{\au{\uep}^{4}}\langle \uep'', A\uep\rangle -
    \frac{|\uep''|^{2}}{\au{\uep)}^{2(\gamma+1)}}\frac{1}{(1+t)^{p}}\right)+
    \nonumber \\
     & & +2(1+t)^{2(p+1)}\left( p\frac{\langle \uep',
     \uep''\rangle}{\au{\uep}^{2(\gamma+1)}}\frac{1}{(1+t)^{p+1}}-
     \frac{\langle A\uep, \uep'\rangle}{\auq{\uep}}
     \frac{\auq{\uep'}}{\auq{\uep}} \right) + \nonumber \\
     & & - 2(\gamma+1)\ep(1+t)^{2(p+1)}\frac{\langle A\uep,
     \uep'\rangle}{\auq{\uep}}
     \frac{|\uep''|^{2}}{\au{\uep}^{2(\gamma+1)}}
     +2(p+1)\frac{R_{\ep}}{1+t};
    \label{derR}
\end{eqnarray}
\begin{eqnarray}
    H_{\ep}' & = & \frac{p+1}{1+t}H_{\ep} -2 (1+t)^{p+1}
    \frac{\auq{\uep'}}{\auq{\uep}}\left(\frac{1}{(1+t)^{p}} + 
    \ep \frac{\langle A\uep, 
     \uep'\rangle}{\auq{\uep}} \right)+
    \nonumber \\
     &  & +2(\gamma - 1)(1+t)^{p+1}\frac{\langle A\uep, 
     \uep'\rangle}{\auq{\uep}} \au{\uep}^{2(\gamma-1)}|A\uep|^{2}.
    \label{derH}
\end{eqnarray}
We are now ready to prove (\ref{SQ}) trough (\ref{SH}).
\subparagraph{\textmd{\textsl{Proof of (\ref{SQ})}}}
Thanks to (\ref{derQ}) we have
$$Q_{\ep}'\leq 
-\frac{1}{\ep}Q_{\ep}\left(\frac{2}{(1+t)^{p}}-\frac{\ep(1-p)}{1+t} + 
2\ep(\gamma+1)\frac{\langle A\uep, 
     \uep'\rangle}{\auq{\uep}}\right) + 
     \frac{2}{\ep}(1+t)^{1-p}\frac{|\uep'||A\uep|}{\auq{\uep}}.$$
 Moreover by (\ref{fond}) and (\ref{1-ep}) we get
 \begin{eqnarray*}
    \frac{2}{(1+t)^{p}}-\frac{\ep(1-p)}{1+t} + 
2\ep(\gamma+1)\frac{\langle A\uep, 
     \uep'\rangle}{\auq{\uep}} & \geq & \frac{1}{(1+t)^{p}}(2 - \ep 
     -2\ep(\gamma+1)K_{0})  \\
     & \geq & \frac{7}{4}\frac{1}{(1+t)^{p}} \geq \frac{1}{(1+t)^{p}}.
\end{eqnarray*}
Hence, using (\ref{fond}) once again, we obtain
\begin{eqnarray*}
   Q_{\ep}'  & \leq &  -\frac{1}{\ep}\frac{1}{(1+t)^{p}}Q_{\ep} + 
   \frac{2}{\ep}\frac{1}{(1+t)^{p}}\sqrt{Q_{\ep}}|A\uep|\au{\uep}^{\gamma - 
   1}(1+t)^{(p+1)/2}  \\
     & \leq &  -\frac{1}{\ep}\frac{1}{(1+t)^{p}}\sqrt{Q_{\ep}} (\sqrt{Q_{\ep}} 
     -2\sqrt{K_{1}}).
\end{eqnarray*}
Therefore, since $Q_{\ep}(0) = Q_{1}(0)$, thesis follows from 
Statement
(\ref{ODEA1}) in Lemma
\ref{lemma:ode}.
\subparagraph{\textmd{\textsl{An intermediate estimate}}}
Thanks to (\ref{SQ}), for all $\alpha(t)> 0$ it holds true that in $[0,S[$:
\begin{eqnarray}
     \frac{|\langle \uep'(t), 
    \uep''(t)\rangle|}{\au{\uep(t)}^{2(\gamma+1)}} & \leq & 
    \frac{1}{2}\alpha(t)  
    \frac{|\uep''(t)|^{2}}{\au{\uep(t)}^{2(\gamma+1)}}  
    + \frac{1}{2\alpha(t)}
    \frac{|\uep'(t)|^{2}}{\au{\uep(t)}^{2(\gamma+1)}}  
    \nonumber  \\
     & \leq & \frac{1}{2}\alpha(t)  
    \frac{|\uep''(t)|^{2}}{\au{\uep(t)}^{2(\gamma+1)}}  + 
    \frac{1}{2\alpha(t)}\frac{L_{2}}{(1+t)^{1-p}}.
    \label{SIn}
\end{eqnarray}
\subparagraph{\textmd{\textsl{Proof of (\ref{SD})}}}
Let us estimate separately the terms in (\ref{derD}).

Using (\ref{fond}) and (\ref{SIn}) with $\alpha(t) = (1+t)(3 + 
2(\gamma+1)K_{0})$ we obtain
\begin{eqnarray*}
    |\varphi_{1}| & \leq & 
    \ep\left[\frac{|\uep''|^{2}}{\au{\uep}^{2(\gamma+1)}}(1+t)^{2p+1} + 
    \frac{|\langle \uep', 
    \uep''\rangle|}{\au{\uep}^{2(\gamma+1)}} (3 + 
    2(\gamma+1)K_{0})(1+t)^{2p}\right]  \\
     & \leq & 
     \ep\left[\frac{|\uep''|^{2}}{\au{\uep}^{2(\gamma+1)}}(1+t)^{2p+1}
     \left(1+ \frac{1}{2}(3+2(\gamma+1)K_{0})^{2}\right) + 
    \frac{1}{2}\frac{L_{2}}{(1+t)^{2-3p}}\right].  
\end{eqnarray*}
Thus from the smallness assumption (\ref{2-ep}) we get
\begin{equation}
    |\varphi_{1}|\leq \frac{1}{16(K_{0}+1)}
    \frac{|\uep''|^{2}}{\au{\uep}^{2(\gamma+1)}}(1+t)^{2p+1}+
    \frac{1}{2}L_{2}(1+t)^{p}.
    \label{Sphi1}
\end{equation}
From (\ref{SIn}) with $\alpha(t) = (1+t)^{p}(8(K_{0}+1))^{-1}$ we have
\begin{equation}
    |\varphi_{3}| \leq \frac{1}{16(K_{0}+1)}
    \frac{|\uep''|^{2}}{\au{\uep}^{2(\gamma+1)}}(1+t)^{2p+1}
    + 4(K_{0}+1)L_{2}(1+t)^{p}.
    \label{Sphi3}
\end{equation}
Moreover from (\ref{SQ}) we get
\begin{equation}
    |\varphi_{4}| \leq \frac{L_{2}}{(1+t)^{1-p}}(1+t)^{p} \leq 
    L_{2}(1+t)^{p}.
    \label{Sphi4}
\end{equation}
Finally replacing (\ref{Sphi1}), (\ref{Sphi3}), 
(\ref{Sphi4}) in  (\ref{derD}), since $\varphi_{2}\leq 0$ we obtain
$$D_{\ep}' \leq \frac{1}{8(K_{0}+1)}
    \frac{|\uep''|^{2}}{\au{\uep}^{2(\gamma+1)}}(1+t)^{2p+1} +
    \left(\frac{11}{2} + 4K_{0}\right)L_{2}(1+t)^{p}.$$
    Hence (\ref{SD}) follows from  a simple integration.
    \subparagraph{\textmd{\textsl{Proof of (\ref{SR})}}}
    Firstly let us estimate some of the terms in (\ref{derR}).
    
    Thanks to (\ref{fond}) we have
    \begin{eqnarray}
        2\gamma \frac{|\langle A\uep, 
     \uep'\rangle|}{\auq{\uep}} \frac{|\langle \uep'', 
     A\uep\rangle|}{\auq{\uep}}& \leq & 2\gamma \frac{K_{0}}{1+t}
     \frac{|\uep''|}{\au{\uep}^{(\gamma+1)}}|A\uep|\au{\uep}^{\gamma-1}
        \nonumber  \\
         & \leq & 
	 \frac{1}{8}
	 \frac{|\uep''|^{2}}{\au{\uep}^{2(\gamma+1)}}\frac{1}{(1+t)^{p}}+
	 \frac{8\gamma^{2}K_{0}^{2}}{(1+t)^{2-p}}|A\uep|^{2}
	 \au{\uep}^{2(\gamma-1)}
        \nonumber\\
         & \leq &  \frac{1}{8}
	 \frac{|\uep''|^{2}}{\au{\uep}^{2(\gamma+1)}}\frac{1}{(1+t)^{p}}+ 
	  \frac{8\gamma^{2}K_{0}^{2}K_{1}}{(1+t)^{3}}.
        \label{1-SR}
    \end{eqnarray}
    Moreover from (\ref{SIn}) with $\alpha(t) = (1+t)/4$ we achieve
    \begin{equation}
        \frac{p}{(1+t)^{p+1}} \frac{|\langle \uep', 
    \uep''\rangle|}{\au{\uep}^{2(\gamma+1)}} \leq 
    \frac{1}{8}
	 \frac{|\uep''|^{2}}{\au{\uep}^{2(\gamma+1)}}\frac{1}{(1+t)^{p}}+
	 \frac{2L_{2}}{(1+t)^{3}}.
        \label{2-SR}
    \end{equation}
     Using once again (\ref{fond}) we have
  \begin{equation}
      \frac{|\langle A\uep, 
     \uep'\rangle|}{\auq{\uep}}\frac{\auq{\uep'}}{\auq{\uep}}
     \leq \frac{K_{0}}{1+t}\frac{\auq{\uep'}}{\auq{\uep}}.
      \label{4-SR}
  \end{equation}
  Finally from (\ref{fond}) and (\ref{1-ep}) we get also
  \begin{eqnarray}
      \ep(\gamma+1)\frac{|\langle A\uep, 
     \uep'\rangle|}{\auq{\uep}} 
     \frac{|\uep''|^{2}}{\au{\uep}^{2(\gamma+1)}} 
     & \leq & \ep \frac{(\gamma + 1)K_{0}}{1+t}
     \frac{|\uep''|^{2}}{\au{\uep}^{2(\gamma+1)}} 
      \nonumber  \\
       & \leq & \frac{1}{8}
	 \frac{|\uep''|^{2}}{\au{\uep}^{2(\gamma+1)}}\frac{1}{(1+t)^{p}}.
      \label{3-SR}
  \end{eqnarray}
 
  Replacing (\ref{1-SR}),  (\ref{2-SR}),  (\ref{4-SR}),  (\ref{3-SR}) 
  in (\ref{derR}) and using (\ref{1-ep}) and (\ref{defsigma1}) we thus obtain
  \begin{eqnarray*}
      R_{\ep}' & \leq & -  
      \frac{|\uep''|^{2}}{\au{\uep}^{2(\gamma+1)}} 
      (1+t)^{p+2}\left(\frac{5}{4}-\frac{2\ep(p+1)}{(1+t)^{1-p}}\right)
      + 
       \\
       & & +2(K_{0}+p+1)\frac{\auq{\uep'}}{\auq{\uep}}(1+t)^{2p+1} +
      \\
       &  & 
       +2(8\gamma^{2}K_{0}^{2}K_{1}+2L_{2})\frac{(1+t)^{p}}{(1+t)^{1-p}}  \\
       & \leq & -\frac{|\uep''|^{2}}{\au{\uep}^{2(\gamma+1)}} 
      (1+t)^{p+2} + 4(K_{0}+1)\frac{\auq{\uep'}}{\auq{\uep}}(1+t)^{2p+1}
     + h_{3}(1+t)^{p}.
  \end{eqnarray*}
Then integrating in $[0,t]$  and using 
(\ref{SD}), since $2p+1 \leq p+2$ we find
\begin{eqnarray}
     & R_{\ep}(t) + \displaystyle \int_{0}^{t}
     \frac{|\uep''(s)|^{2}}{\au{\uep(s)}^{2(\gamma+1)}} 
      (1+s)^{p+2}ds 
      &   \nonumber\\
     & \leq 4(K_{0}+1)
       \displaystyle \int_{0}^{t}\frac{\auq{\uep'(s)}}{\auq{\uep(s)}}(1+s)^{2p+1} ds 
      + \frac{h_{3}}{1+p}(1+t)^{p+1} + R_{\ep}(0) &   \nonumber\\
     & \leq  -4(K_{0}+1)\left[\ep \displaystyle\frac{\langle \uep'(t), 
    \uep''(t)\rangle}{\au{\uep(t)}^{2(\gamma+1)}}(1+t)^{2p+1} - D_{\ep}(0)\right] + 
  & \nonumber\\
   & +\displaystyle  \frac{1}{2}\displaystyle \int_{0}^{t}
     \frac{|\uep''(s)|^{2}}{\au{\uep(s)}^{2(\gamma+1)}} 
      (1+s)^{p+2}ds+ &  \nonumber \\
     &  \displaystyle  +  8(K_{0}+1)(3+2K_{0})L_{2}(1+t)^{p+1}+  
     \frac{h_{3}}{1+p}(1+t)^{p+1} + R_{\ep}(0). \nonumber \label{stimaR}
\end{eqnarray}
From this inequality and (\ref{defsigma1})  it follows that:
\begin{eqnarray}
    & R_{\ep}(t) + \displaystyle \frac{1}{2}\int_{0}^{t}
   \frac{|\uep''(s)|^{2}}{\au{\uep(s)}^{2(\gamma+1)}} 
      (1+s)^{p+2}ds   \leq  R_{\ep}(0) + 4(K_{0}+1)|D_{\ep}(0)| +&
    \nonumber \\
     &    +4(K_{0}+1)\ep \displaystyle\frac{|\langle \uep'(t), 
    \uep''(t)\rangle|}{\au{\uep(t)}^{2(\gamma+1)}}(1+t)^{2p+1} + 
    h_{4}(1+t)^{p+1}.&
    \label{5-SR}
\end{eqnarray}
Let us now estimate the terms in the right hand side. Let us remark that from (\ref{SIn}) with $\alpha(t) = 
4\ep(K_{0}+1)(1+t)$ we get
$$
    4(K_{0}+1)\ep \displaystyle\frac{|\langle \uep', 
    \uep''\rangle|}{\au{\uep}^{2(\gamma+1)}}(1+t)^{2p+1} \leq 
    8\ep^{2}(K_{0}+1)^{2}\frac{|\uep''|^{2}(1+t)^{2p+2}}{\au{\uep}^{2(\gamma+1)}} 
       + \frac{L_{2}}{2}(1+t)^{3p-1}.
    $$
      Moreover $3p - 1 \leq 2p\leq p+1$, hence using also 
      (\ref{1-ep}) we 
      deduce 
       \begin{equation}
           4(K_{0}+1)\ep \displaystyle\frac{|\langle \uep', 
    \uep''\rangle|}{\au{\uep}^{2(\gamma+1)}}(1+t)^{2p+1} \leq  
    \frac{1}{2}R_{\ep} +  \frac{1}{2}L_{2}(1+t)^{p+1}.
           \label{6-SR}
       \end{equation}
Furthermore we have
\begin{eqnarray}
    |D_{\ep}(0)| & \leq & \ep
    |u''(0)|\frac{|u_{1}|}{\au{u_{0}}^{2(\gamma+1)}} =
    |u_{1}+\au{u_{0}}^{2\gamma}Au_{0}|\frac{|u_{1}|}{\au{u_{0}}^{2(\gamma+1)}}
    \nonumber \\
    & \leq & (|u_{1}| +\au{u_{0}}^{2\gamma}|Au_{0}|)
    \frac{|u_{1}|}{\au{u_{0}}^{2(\gamma+1)}}.
    \label{7-SR}
\end{eqnarray}
Plugging (\ref{6-SR}) and  (\ref{7-SR}) in  (\ref{5-SR}) and 
recalling the definition in
(\ref{defL2}) we get
\begin{eqnarray*}
   &  \displaystyle\frac{1}{2} R_{\ep}(t) + \displaystyle \frac{1}{2}\int_{0}^{t}
   \frac{|\uep''(s)|^{2}}{\au{\uep(s)}^{2(\gamma+1)}} 
      (1+s)^{p+2}ds  & \\
      &\leq  \left( R_{\ep}(0) + 
      4(K_{0}+1)|D_{\ep}(0)| + h_{4} + 
      \displaystyle \frac{1}{2}L_{2}\right)(1+t)^{p+1} & \\
     & \leq \left( R_{\ep}(0) + \displaystyle \frac{1}{2}L_{3}\right)(1+t)^{p+1}.   &
\end{eqnarray*}
Therefore inequality (\ref{SR}) is proved.
\subparagraph{\textmd{\textsl{Proof of (\ref{SH})}}} If $\gamma = 1$, thesis, eventually 
for smaller values of $\ep$, follows from Theorem 2.2 in \cite{jde2} 
(in particular it is a consequence of (3.52)). Then from now let us 
assume that $\gamma > 1$. 
Let us set
$$F_{\ep}(t):= \ep \frac{\auq{\uep'(t)}}{\auq{\uep(t)}} (1+t)^{p+1}, 
\hspace{1em}G_{\ep}(t):=
    \au{\uep(t)}^{2(\gamma-1)}|A\uep(t)|^{2} (1+t)^{p+1}$$
    so that
    $H_{\ep} = F_{\ep} + G_{\ep}$.
    Hence exploiting (\ref{fond}) in (\ref{derH})  we obtain:
    $$H_{\ep}' \leq -\frac{2}{\ep}F_{\ep}\left(\frac{1}{(1+t)^{p}} 
    -\ep(K_{0}+p+1)\frac{1}{1+t} 
   \right) +\frac{p+1}{1+t}G_{\ep} + 2(\gamma -1) \frac{\langle A\uep, 
     \uep'\rangle}{\auq{\uep}} G_{\ep}.$$
     Thanks to (\ref{1-ep}) we have $1 - \ep(K_{0}+2) \geq 1/2$, 
     therefore we get
   \begin{equation}
       H_{\ep}' \leq -\frac{F_{\ep}}{\ep}\frac{1}{1+t} 
  +\frac{2}{1+t}G_{\ep} + 2(\gamma -1) \frac{\langle A\uep, 
     \uep'\rangle}{\auq{\uep}} G_{\ep}.
       \label{1-SH}
   \end{equation}  
   Now let us recall that
   $$\uep' = -(1+t)^{p}(\ep\uep'' + \au{\uep}^{2\gamma}A\uep),$$
   hence
   $$
       \frac{\langle A\uep, 
     \uep'\rangle}{\auq{\uep}} = -\ep\frac{\langle A\uep, 
     \uep''\rangle}{\auq{\uep}}(1+t)^{p} - \frac{1}{1+t}G_{\ep}.
       $$
   Plugging this identity in (\ref{1-SH}), we arrive at
    \begin{equation}
        H_{\ep}'\leq -\frac{1}{1+t}\left(\frac{F_{\ep}}{\ep} 
-2G_{\ep} + 2(\gamma - 1)G_{\ep}^{2} +2\ep(\gamma - 
	1)(1+t)^{p+1}G_{\ep}\frac{\langle A\uep, 
     \uep''\rangle}{\auq{\uep}}\right).
        \label{3-SH}
    \end{equation}
   Let us estimate the last term in (\ref{3-SH}).  Using (\ref{SR}) we  obtain 
    \begin{eqnarray}
        2\ep(\gamma - 
	1)(1+t)^{p+1}\frac{| A\uep| 
    |\uep''|}{\auq{\uep}} & = & 2(\gamma - 1)\ep 
    \frac{|\uep''|}{\au{\uep}^{\gamma+1}}(1+t)^{(1+p)/2}\sqrt{G_{\ep}}
        \nonumber  \\
         & \leq &2(\gamma - 1) \sqrt{\ep}\sqrt{G_{\ep}}\sqrt{2R_{\ep}(0) + 
	 L_{3}}
        \nonumber  \\
         & \leq & 2(\gamma - 
	 1)\sqrt{G_{\ep}}\left(\sqrt{2\ep 
	 R_{\ep}(0) }+
	 \sqrt{\ep}\sqrt{L_{3}}\right).
        \label{4-SH}
    \end{eqnarray}
    By  the definition of $h_{1}$ moreover it follows that:
    \begin{eqnarray}
        2\ep 
	 R_{\ep}(0)  & = 
	 &2\ep^{2}\frac{|\uep''(0)|^{2}}
	 {\au{u_{0}}^{2(\gamma+1)}}+2\ep\frac{\auq{u_{1}}}{\auq{u_{0}}}
        \nonumber  \\
	 & \leq &
	4 (|u_{1}|^{2}+\au{u_{0}}^{4\gamma}|Au_{0}|^{2})
	 \frac{1}{\au{u_{0}}^{2(\gamma+1)}}
	 +2\ep\frac{\auq{u_{1}}}{\auq{u_{0}}}
	\nonumber  \\
         & = & h_{1}  +2\ep\frac{\auq{u_{1}}}{\auq{u_{0}}}.
         \label{5-SH}
    \end{eqnarray}
    Using (\ref{3-ep})  and (\ref{5-SH}), from (\ref{4-SH}) we get
    $$
        2\ep(\gamma - 1)(1+t)^{p+1}\frac{|\langle A\uep, 
    \uep''\rangle|}{\auq{\uep}}  \leq 2(\gamma - 1)\sqrt{G_{\ep}} 
    (\sqrt{h_{1}} + 1).
       $$
    Plugging this estimate in (\ref{3-SH}), since $\ep \leq 1$ we 
    finally achieve
    $$H_{\ep}'  \leq -\frac{1}{1+t}\left(F_{\ep} + 2(\gamma - 1)G_{\ep}^{2} 
	-2G_{\ep} 
	-2(\gamma-1)G_{\ep}^{3/2}(\sqrt{h_{1}} + 1)\right).$$
Since $H_{\ep}(0)\leq H_{1}(0)$, inequality (\ref{SH}) follows recalling the definition of 
$L_{1}$ and applying Lemma \ref{lemma:ode-GF} 
with 
$$a=2(\gamma-1), \hspace{1em} b = 2, \hspace{1em} 
c= 2(\gamma -1)(\sqrt{h_{1}} + 1).$$ \qed
\subsubsection{Existence of global solutions}

\subparagraph{\textmd{\textsl{Local maximal solutions}}}

Problem (\ref{pbm:h-eq}), (\ref{pbm:h-data}) admits a unique local-in-time
solution, and this solution can be continued to a solution defined in 
a maximal interval $[0,T[$, where either $T=+\infty$, or
\begin{equation}
	\limsup_{t\to T^{-}}\left(|A^{1/2}\uep'(t)|^{2}+
	|A\uep(t)|^{2}\right)=+\infty,
	\label{limsup}
\end{equation}
or
\begin{equation}
	\liminf_{t\to T^{-}}|A^{1/2}\uep(t)|^{2}=0.
	\label{liminf}
\end{equation}

We omit the proof of this standard result.  The interested reader is
referred to \cite{gg:k-dissipative}.
\subparagraph{\textmd{\textsl{The standard conserved energy}}}
Let us recall that problem (\ref{pbm:h-eq}), (\ref{pbm:h-data}) 
admits a first order conserved energy, that is
\begin{equation}
    \ep|\uep'(t)|^{2} + \frac{\au{\uep(t)}^{2(\gamma+1)}}{\gamma+1} 
    +2\int_{0}^{t}\frac{|\uep'(s)|^{2}}{(1+s)^{p}}ds = \ep|u_{1}|^{2}
    +\frac{\au{u_{0}}^{2(\gamma+1)}}{\gamma+1}. 
    \label{conserv}
\end{equation}
Therefore $\auq{\uep(t)}$ is bounded independently from $\ep\leq 1$.
\subparagraph{\textmd{\textsl{Global solutions}}}
We want to apply  Proposition \ref{priori}.  To this end let us set:
\begin{equation}
    K_{1}:= L_{1}+1+H_{1}(0), \hspace{1em} K_{0}:=\sqrt{\max\{4K_{1}, 
Q_{1}(0)\}K_{1}} +   \frac{|\langle Au_{0}, 
     u_{1}\rangle|}{\auq{u_{0}}}.
    \label{defk}
\end{equation}
For such choices of $K_{0}$ and $K_{1}$ let us define
$$S:=\sup\{\tau\in[0,T[:\; \mbox{(\ref{fond1}), (\ref{fond}) are 
verified for all $ t\in [0,\tau]$}\}. $$

Firstly let us remark that since $\auq{u_{0}} > 0$, for our choices of $K_{0}$
and $K_{1}$  and $\ep \leq 1$ we 
have $S > 0$. From now furthermore we assume that $\ep$ verifies the smallness 
conditions of Proposition \ref{priori}. Thus in $[0,S[$ Proposition \ref{priori} 
holds true.   

We want to prove that $S = T$. 

Let us assume by 
contradiction that $S < T$. Then by the regularity properties of $\uep$  all the 
estimates in Proposition \ref{priori} hold true in $[0,S]$. Moreover at least one of the following is 
verified:
\begin{equation}
    \auq{\uep(S)} = 0,
    \label{exc1}
\end{equation}
\begin{equation}
    \frac{|\langle \uep'(S), 
    A\uep(S)\rangle|}{\auq{\uep(S)}} =
   \frac{K_{0}}{1+S}, 
    \label{exc2}
\end{equation}
\begin{equation}
    \au{\uep(S)}^{2(\gamma-1)}|A\uep(S)|^{2} =
   \frac{K_{1}}{(1+S)^{p+1}}.
    \label{exc3}
\end{equation}

\subparagraph{\textmd{\em{Equality (\ref{exc1}) is false}}}
Let us set $y(t):=\auq{\uep(t)}$. Hence by (\ref{fond}) in $[0,S[$ we 
have:
$$
    \frac{y'(t)}{y(t)}\geq -\frac{2K_{0}}{1+t}
   $$
 therefore 
   \begin{equation}
       y(t) \geq y(0) e^{-2K_{0}\log(1+t)} = 
       \frac{\auq{u_{0}}}{(1+t)^{2K_{0}}};
       \label{decy}
   \end{equation}
   in particular  $\auq{\uep(S)} = y(S)  > 0$. 
\subparagraph{\textmd{\em{Equality (\ref{exc2}) is false}}}
From (\ref{SQ}) and (\ref{SH}), recalling (\ref{defk}) and the definition of
$L_{2}$  in Proposition \ref{priori} we indeed have
\begin{eqnarray*}
     \frac{|\langle \uep'(S), 
    A\uep(S)\rangle|}{\auq{\uep(S)}} &  \leq & \frac{|\uep'(S)|}{\au{\uep(S)}^{\gamma+1}}
    |A\uep(S)| \au{\uep(S)}^{\gamma-1}  \\
     & \leq & \frac{\sqrt{L_{2}}}{(1+S)^{(1-p)/2}} 
     \frac{\sqrt{L_{1}}}{(1+S)^{(1+p)/2}}
     =  
     \frac{\sqrt{L_{2}L_{1}}}{1+S}<\frac{\sqrt{L_{2}K_{1}}}{1+S} 
      \leq  \frac{K_{0}}{1+S}.
\end{eqnarray*}
\subparagraph{\textmd{\em {Equality (\ref{exc3}) is false}}} This is an immediate 
consequence of (\ref{SH}).
\subparagraph{\textmd{\textsl{Conclusion}}}
We have proved that if $\ep$ is small enough hence $S = T$ and in 
$[0,T[$ Proposition \ref{priori} holds true. We need only to prove 
that $T = +\infty$. If it is not the case hence (\ref{limsup}) or 
(\ref{liminf}) hold true. Nevertheless (\ref{liminf})  is excluded 
by (\ref{decy}). Now let us roll out  (\ref{limsup}). From (\ref{conserv}) 
we know that $\auq{\uep}$ is bounded from above, hence from (\ref{SH}) 
 we deduce that $\auq{\uep'}$ is bounded. Moreover thanks to 
 (\ref{decy})  $\auq{\uep}$ is bounded also from below, thus using once 
 again   (\ref{SH})  we get that  $|A\uep|^{2}$ is bounded. Hence (\ref{limsup})
 is false.

\subsubsection{Decay estimates}
Since now  inequalities  (\ref{SH}) and  (\ref{SQ})
hold true in $[0,+\infty[$, we have already proved (\ref{D2}) and (\ref{D3}).
The estimate from below in (\ref{D1}) is a consequence of 
Proposition 3.3 (see (3.21)) in \cite{jde2}, because we have proved 
that 
$$\frac{|\langle \uep'(t), 
    A\uep(t)\rangle|}{\auq{\uep(t)}} \leq \frac{K_{0}}{1+t} 
    \hspace{1em} \forall t \geq 0.$$
   \subparagraph{\textmd{\textsl{Proof of  (\ref{D0}), (\ref{D4}) and estimate from above 
    in (\ref{D1})}}} Let us recall that we have already supposed that the 
    smallness assumption (\ref{1-ep}) is verified.
   We work as in \cite{jde2}, Section 3.4, hence we skip the details.
   Let us set
   $$\mathcal{D}_{\ep}(t):= 
   \ep(1+t)^{p}\langle\uep'(t),\uep(t)\rangle + \frac{1}{2}\left(1 - 
  \frac{ \ep p}{(1+t)^{1-p}}\right)|\uep(t)|^{2}. $$
  Since 
   $$\mathcal{D}_{\ep}' = -(1+t)^{p}\au{\uep}^{2(\gamma+1)} +
   \ep(1+t)^{p}|\uep'|^{2}  + 
   \ep\frac{p(1-p)}{2}\frac{|\uep|^{2}}{(1+t)^{2-p}},$$
   a simple integration gives
   \begin{eqnarray*}
     \int_{0}^{t}(1+s)^{p}\au{\uep(s)}^{2(\gamma+1)}ds    & = &
     \mathcal{D}_{\ep}(0) - \mathcal{D}_{\ep}(t) +
   \ep \int_{0}^{t}(1+s)^{p}|\uep'(s)|^{2}ds    \\
        &  &  +
   \ep\frac{p(1-p)}{2}\int_{0}^{t}\frac{|\uep(s)|^{2}}{(1+s)^{2-p}}ds.
\end{eqnarray*}
  Moreover, since $1-3\ep\geq 1/4$, it holds true that
  $$  - \mathcal{D}_{\ep}(t) \leq 
  \frac{1}{4}\ep(1+t)^{p+1}|\uep'(t)|^{2} - 
  \frac{1}{8}|\uep(t)|^{2},$$
  hence we get:
  \begin{eqnarray}
     & \displaystyle \frac{1}{8}|\uep(t)|^{2} +  
      \int_{0}^{t}(1+s)^{p}\au{\uep(s)}^{2(\gamma+1)}ds   \leq  
      |\mathcal{D}_{\ep}(0)| + 
      \frac{1}{4}\ep(1+t)^{p+1}|\uep'(t)|^{2}+&
      \nonumber  \\
       &  + \displaystyle \ep \int_{0}^{t}(1+s)|\uep'(s)|^{2}ds +
       \ep\frac{p(1-p)}{2}\int_{0}^{t}\frac{|\uep(s)|^{2}}{(1+s)^{2-p}}ds.& 
      \label{STD}
  \end{eqnarray}
  Let us now define
      $$\mathcal{E}_{\ep}(t):=  \left(\ep|\uep'(t)|^{2} + 
    \frac{\au{\uep(t)}^{2(\gamma+1)}}{\gamma+1}\right)(1+t)^{p+1}.$$  
A simple computation gives
$$\mathcal{E}_{\ep}' = -(1+t)\left(2 - 
\frac{\ep(p+1)}{(1+t)^{1-p}} \right)|\uep'|^{2}+ 
\frac{p+1}{\gamma+1}(1+t)^{p}\au{\uep}^{2(\gamma+1)}.$$
Integrating in $[0,t]$ and using (\ref{STD}) we arrive at
\begin{eqnarray*}
  &  \displaystyle (1+t)^{p+1}\left(1 - 
  \frac{p+1}{4(\gamma+1)}\right)\ep |\uep'(t)|^{2} + 
    \frac{\au{\uep(t)}^{2(\gamma+1)}}{\gamma+1}(1+t)^{p+1} \leq
      &   \\
     &\displaystyle  \leq  \mathcal{E}_{\ep}(0) -
     \left(2 - \ep (p+1) - \ep\frac{p+1}{\gamma+1}\right)
     \int_{0}^{t}(1+s)|\uep'(s)|^{2}ds +  &\\
     & \displaystyle  + \frac{p+1}{\gamma+1}\left( |\mathcal{D}_{\ep}(0)| -
     \displaystyle \frac{1}{8}|\uep(t)|^{2}
    + \ep\frac{p(1-p)}{2}\int_{0}^{t}
    \frac{|\uep(s)|^{2}}{(1+s)^{2-p}}ds \right). &
\end{eqnarray*}
Since $2 -2\ep(1+p)\geq 1$, then it holds true that
\begin{eqnarray}
   &\displaystyle  \frac{1}{2}\mathcal{E}_{\ep}(t) + \int_{0}^{t}(1+s)|\uep'(s)|^{2}ds 
    +  \frac{1}{8} \frac{p+1}{\gamma+1}|\uep(t)|^{2}  
     \leq & 
    \nonumber \\
     &\displaystyle \mathcal{E}_{1}(0)+2|\mathcal{D}_{\ep}(0)| + 
     \ep\frac{p+1}{\gamma+1}\frac{p(1-p)}{2}\int_{0}^{t}
    \frac{|\uep(s)|^{2}}{(1+s)^{2-p}}ds. & 
    \label{1-SE}
\end{eqnarray}
In particular we have
$$|\uep(t)|^{2}  \leq 
8\frac{\gamma+1}{p+1}(\mathcal{E}_{1}(0)+2|\mathcal{D}_{\ep}(0)|)
+ 4\ep(1-p)\int_{0}^{t}
    \frac{|\uep(s)|^{2}}{(1+s)^{2-p}}ds.$$
    From the Gronwall's Lemma we hence get
    $$|\uep(t)|^{2}  \leq 16\frac{\gamma+1}{p+1}
    (\mathcal{E}_{1}(0)+2|\mathcal{D}_{\ep}(0)|),$$
    and finally
    $$(1-p)\int_{0}^{t}
    \frac{|\uep(s)|^{2}}{(1+s)^{2-p}}ds \leq 16\frac{\gamma+1}{p+1}
    (\mathcal{E}_{1}(0)+2|\mathcal{D}_{\ep}(0)|).$$
   Now we go back to (\ref{1-SE}) and from the previous inequality we obtain
   $$\frac{1}{2}\mathcal{E}_{\ep}(t) + \int_{0}^{t}(1+s)|\uep'(s)|^{2}ds 
    +  \frac{1}{8} \frac{p+1}{\gamma+1}|\uep(t)|^{2}  \leq 
    9(\mathcal{E}_{1}(0) + |u_{1}||u_{0}| + |u_{0}|^{2}).$$
    By this inequality all the  estimates we 
    look for immediately follow.
    \qed
    \subsection{Proof of Theorem \ref{thm:dg-error}}\label{error}
Also in this proof in most cases we omit the dependence of $u$, 
$\uep$, $\roep$, $\rep$ and $\tetep$ from t.
Moreover from now on we assume that $\ep$ verifies the smallness assumptions of 
 Theorem \ref{th:ex} in such a way that $\uep$ is globally well 
 defined. 
 
  Let us recall that $\rep$ and $\roep$  verify the 
   following problems:
   \begin{equation}
       \left\{
       \begin{array}{l}
       \displaystyle \ep \rep'' + \au{u}^{2\gamma}A\roep 
       +\frac{1}{(1+t)^{p}}\rep'= -\ep u''+ (\au{u}^{2\gamma} - 
       \au{\uep}^{2\gamma})A\uep      \\
           \rep(0) = \rep'(0) =0;
       \end{array}\right.
       \label{cpr}
   \end{equation}
   and 
   \begin{equation}
      \left\{
       \begin{array}{l}
       \displaystyle \ep \roep'' + (\au{\uep}^{2\gamma}A\uep - \au{u}^{2\gamma}Au) 
       +\frac{1}{(1+t)^{p}}\roep'= -\ep u'' \\ \roep(0) = 
       0,\;\roep'(0) = w_{0};
         \end{array}\right.
       \label{cpro}
   \end{equation}
   where $w_{0}$, $\uep$ and $u$ are defined in (\ref{theta-data}), 
   (\ref{pbm:h-eq}) and (\ref{pbm:par}) respectively.
    \subsubsection{Fundamental energies}
   Below we define the energies  we use in the proof of Theorem \ref{thm:dg-error}.
   
   Let us set
    \begin{eqnarray}
       D_{\rho}(t)&:= &\int_{0}^{t}\langle \au{\uep(s)}^{2\gamma}A\uep(s) - 
       \au{u(s)}^{2\gamma}Au(s), \roep(s)\rangle (1+s)^{p} ds+ \nonumber \\
       & & +\ep\langle\roep'(t),\roep(t)\rangle (1+t)^{p}+ 
       \frac{1}{2}|\roep(t)|^{2}(1 - \ep p(1+t)^{p-1}) ;   
       \label{Dro}
   \end{eqnarray}
   \begin{equation}
       E_{\rho}(t):= (\ep|\rep'(t)|^{2}+ 
       \au{u(t)}^{2\gamma}\auq{\roep(t)})(1+t)^{2p};
       \label{Ero}
   \end{equation}
  \begin{equation}
       F_{\rho}(t):= \ep\frac{|\rep'(t)|^{2}}{\au{u(t)}^{2\gamma}}+ 
       \auq{\roep(t)}.
       \label{Fro}
   \end{equation}
   In the proposition below we recollect all the inequalities verified by 
   these energies we need.
    \begin{prop} [Basic error estimates]\label{energyest}
Let $A$ be a non negative operator and let us assume that $0\leq p 
\leq 1$, $\gamma \geq 1$. Moreover let us suppose that  $(u_{0},u_{1})\in 
D(A^{3/2})\times D(A^{1/2})$. Then there exists $\ep_{0}>0$ such that for all $\ep\leq \ep_{0}$  
	and for all $t\geq 0$	we have
	\begin{eqnarray}&
	   \displaystyle  |\roep(t)|^{2}+\int_{0}^{t}(\au{u(s)}^{2\gamma}+\au{\uep(s)}^{2\gamma})
	    \auq{\roep(s)}(1+s)^{p}ds \leq & \nonumber \\
	    &\displaystyle \leq
	    \gamma_{6}\ep^{2}+8\ep\int_{0}^{t}|\roep'(s)|^{2}(1+s)^{p}ds;&
	    \label{stDro1}
	\end{eqnarray}
	\begin{equation}
	    F_{\rho}(t) + 
	    \int_{0}^{t}\frac{|\rep'(s)|^{2}}{\au{u(s)}^{2\gamma}}
	    \frac{1}{(1+s)^{p}} ds
	    \leq \gamma_{7}\ep^{2}+\gamma_{8}
	    \ep\int_{0}^{t}|\roep'(s)|^{2}(1+s)^{p}ds;
	    \label{stFro1}
	\end{equation}
       if moreover (\ref{H1cip}) is verified then
       \begin{equation}
           E_{\rho}(t) + \int_{0}^{t}|\rep'(s)|^{2}(1+s)^{p}ds \leq 
	   \gamma_{9}\ep;
           \label{stEro1}
       \end{equation}
       where all  constants do not depend on $\ep$ and $t$ but 
       only on the initial data.
    \end{prop}
    \prf To begin with, we compute the time's derivatives of 
    (\ref{Ero}), (\ref{Fro})
     and (\ref{Dro}).
   Using (\ref{cpr}) it is easy to see that
   \begin{eqnarray}
       E_{\rho}' & = & -|\rep'|^{2}(1+t)^{p}\left(2 - \frac{2\ep 
       p}{(1+t)^{1-p}}\right) + 
       2p(1+t)^{2p-1}\au{u}^{2\gamma}\auq{\roep} 
       +
       \nonumber  \\
        &  & +2(1+t)^{2p}\au{u}^{2\gamma}\langle A\roep, \tetep'\rangle 
	-2\ep\langle u'',\rep'\rangle (1+t)^{2p} +
       \nonumber  \\
        &  & 
 +2(1+t)^{2p}(\au{u}^{2\gamma}- \au{\uep}^{2\gamma})\langle
A\uep,\rep'\rangle +  \nonumber  \\
& &
-2\gamma(1+t)^{3p}\au{u}^{4\gamma-2}|Au|^{2}\auq{\roep} \nonumber
\\
        & =: & S_{1}+S_{2}+S_{3}+S_{4}+S_{5}+S_{6},
       \label{DEro}
   \end{eqnarray}
   and
   \begin{eqnarray}
       F_{\rho}' & = & 
       -\frac{|\rep'|^{2}}{\au{u}^{2\gamma}}\left(\frac{2}{(1+t)^{p}} 
       -2\gamma\ep\au{u}^{2(\gamma-1)}|Au|^{2}(1+t)^{p}\right) 
       + 2 \langle A\roep, \tetep'\rangle +
       \nonumber\\
        &  & +2\frac{\au{u}^{2\gamma} - 
	\au{\uep}^{2\gamma}}{\au{u}^{2\gamma}}\langle A\uep,\rep'\rangle
	-2\ep\langle u'',\rep'\rangle\frac{1}{\au{u}^{2\gamma}}. 
       \label{DFro}  
   \end{eqnarray}
   Conversely using (\ref{cpro}) we have
   \begin{equation}
       D_{\rho}'= \ep|\roep'|^{2}(1+t)^{p} -\ep\langle 
       u'',\roep\rangle(1+t)^{p} + \frac{p}{2}(1-p)\ep(1+t)^{p-2}|\roep|^{2}.
       \label{DDro}
   \end{equation}
  Moreover from now on let us assume that $\ep_{0}$ verifies also these 
  assumptions (recall that  we have already supposed that 
  $\ep_{0}$ satisfies the smallness conditions in Theorem \ref{th:ex})
  \begin{equation}
      (1+p)\ep_{0} \leq \frac{1}{4},\hspace{1em} 2\gamma 
      \gamma_{4}\ep_{0} \leq \frac{1}{4}
      \label{hep1}
  \end{equation}
  where $\gamma_{4}$ is the constant in (\ref{dv4}).
  In the following we denote by $c_{i}$ various constants that depend 
  only on the initial data.
  Moreover let us set
  $$\phi_{\rho}(t):= (\au{u}^{2\gamma}+
     \au{\uep}^{2\gamma})\auq{\roep}.$$
  \subparagraph{\textmd{\textsl{Preliminary estimates}}}
  Thanks to  Lagrange's Theorem for all $t\geq 0$ there exists 
  $\xi_{t}$ in the interval with end points
  $\auq{u(t)}$ and $\auq{\uep(t)}$ such that 
 \begin{eqnarray*}
     \au{u(t)}^{2\gamma}- \au{\uep(t)}^{2\gamma} & =
     &   \gamma \xi_{t}^{\gamma-1}(\auq{u(t)}-\auq{\uep(t)})  \\
      & = &  -\gamma  
  \xi_{t}^{\gamma-1}\langle A^{1/2}(u(t)+\uep(t)),A^{1/2}\roep(t)\rangle.
 \end{eqnarray*}
Since it is clear that
$$\xi_{t}^{\gamma-1} 
\leq\au{u(t)}^{2(\gamma-1)}+\au{\uep(t)}^{2(\gamma-1)},$$
then
\begin{eqnarray}
    &( \au{u}^{2\gamma}- \au{\uep}^{2\gamma})^{2}  \leq 
    \gamma^{2}\xi_{t}^{2(\gamma-1)}\auq{(u+\uep)}\auq{\roep}& 
    \nonumber  \\
     &  \leq 
     2\gamma^{2}(\au{u}^{2(\gamma-1)}+
     \au{\uep}^{2(\gamma-1)})^{2}(\auq{u}+\auq{\uep})\auq{\roep}&
    \nonumber  \\
     & \leq 6\gamma^{2} (\au{u}^{2(\gamma-1)}+
     \au{\uep}^{2(\gamma-1)}) \phi_{\rho}.& 
    \label{sdiff}  
\end{eqnarray}
Moreover computing the time's derivative of (\ref{pbm:par}) we get:
\begin{equation}
    u''= -p(1+t)^{p-1}\au{u}^{2\gamma}Au + 
    (1+t)^{2p}(\au{u}^{4\gamma}A^{2}u+ 2\gamma
    \au{u}^{4\gamma - 2}|Au|^{2}Au),
    \label{eqv''}
\end{equation}
thus 
\begin{eqnarray}
    |u''|^{2} & \leq  & 3(1+t)^{2(p-1)}\au{u}^{4\gamma}|Au|^{2} + 
    3(1+t)^{4p}\au{u}^{8\gamma}|A^{2}u|^{2} \nonumber \\
     &  & +12\gamma^{2}(1+t)^{4p}\au{u}^{8\gamma-4}|Au|^{6}.
     \label{v''-1}
\end{eqnarray}
From (\ref{dv4}) we also deduce that
\begin{equation}
    \au{u}^{4(\gamma-1)}|Au|^{4} \leq 
    \frac{\gamma_{4}^{2}}{(1+t)^{2p+2}}.
    \label{stvro}
\end{equation}
We can now estimate the last term in (\ref{v''-1}) using 
(\ref{stvro}), so finally we get
\begin{equation}
    |u''|^{2} \leq c_{1}(1+t)^{2p-2}\au{u}^{4\gamma}|Au|^{2} 
    + 3(1+t)^{4p}\au{u}^{8\gamma}|A^{2}u|^{2}.
    \label{sv''}
\end{equation}
Now we are ready  to prove (\ref{stDro1}), (\ref{stFro1}), (\ref{stEro1}).
\subparagraph{\textmd{\textsl{Proof of (\ref{stDro1})}}}
Thanks to Lemma \ref{mprop} with $m(r) = r^{\gamma}$ we have
$$ \frac{1}{2}\phi_{\rho}\leq \langle \au{\uep}^{2\gamma}A\uep - 
       \au{u}^{2\gamma}Au, \roep \rangle
$$
hence integrating (\ref{DDro}) in $[0,t]$  we get
\begin{eqnarray}
   & \displaystyle \frac{1}{2}
   \int_{0}^{t}\phi_{\rho}(s)(1+s)^{p}ds   + 
       \frac{1}{2}|\roep(t)|^{2}(1 - \ep p(1+t)^{p-1})& 
    \nonumber  \\
     & \displaystyle\leq - 
       \ep\langle\roep'(t),\roep(t)\rangle (1+t)^{p}
       +\ep\int_{0}^{t}|\roep'(s)|^{2}(1+s)^{p} ds  +& 
    \nonumber  \\
     & \displaystyle-\ep\int_{0}^{t}\langle 
       u''(s),\roep(s)\rangle(1+s)^{p}ds + \frac{p}{2}(1-p)\ep\int_{0}^{t}(1+s)^{p-2}|\roep(s)|^{2}  ds& 
    \nonumber \\
     &\displaystyle  =:\psi_{1}+\psi_{2}+\psi_{3}+\psi_{4}. & 
    \label{stDro1-1}
\end{eqnarray}
Let us now estimate $\psi_{1}$, $\psi_{3}$, $\psi_{4}$.

From (\ref{D1}) and (\ref{D3}) we obtain
$$|\uep'|^{2} \leq \frac{C_{2}}{(1+t)^{1-p}}
\frac{C_{2}^{\gamma+1}}{(1+t)^{1+p}} = 
\frac{C_{2}^{\gamma+2}}{(1+t)^{2}},$$
and from (\ref{dv4}), (\ref{dv2})  we have
$$|u'|^{2}= (1+t)^{2p}\au{u}^{4\gamma}|Au|^{2} \leq 
\gamma_{4}(1+t)^{2p}\frac{\au{u}^{2(\gamma+1)}}{(1+t)^{p+1}}\leq 
\frac{c_{2}}{(1+t)^{2}}.$$
    Therefore, recalling that $|\roep'|^{2}\leq 2 
    (|\uep'|^{2}+|u'|^{2})$ we get 
    \begin{equation}
        |\psi_{1}|\leq \ep^{2}|\roep'|^{2}(1+t)^{2p} + 
	\frac{1}{4}|\roep|^{2} \leq c_{3}\ep^{2}+ \frac{1}{4}|\roep|^{2}. 
        \label{spsi1ro}
    \end{equation}
    Now let us estimate $\psi_{3}$. From (\ref{eqv''}) we deduce that
   \begin{eqnarray*}
       \langle u'',\roep\rangle   &  = & 
       -\frac{p}{(1+t)^{1-p}
    }\au{u}^{2\gamma}\langle A^{1/2}u, A^{1/2}\roep\rangle  + (1+t)^{2p} \au{u}^{4\gamma}\langle
    A^{3/2}u, A^{1/2}\roep\rangle +  \\
        &  & +2\gamma(1+t)^{2p}\au{u}^{4\gamma-2}|Au|^{2}
	\langle A^{1/2}u,A^{1/2}\roep\rangle,
   \end{eqnarray*}
  hence
  \begin{eqnarray*}
      \ep|\langle u'',\roep\rangle |(1+t)^{p} & \leq & 
      \frac{1}{4}\au{u}^{2\gamma}\auq{\roep}(1+t)^{p} + 
      3\ep^{2}\au{u}^{2(\gamma+1)} (1+t)^{3p-2} +  \\
       &  & +3\ep^{2}(1+t)^{5p}\au{u}^{6\gamma}|A^{3/2}u|^{2}+  \\
       &  & +12\gamma^{2}\ep^{2}(1+t)^{5p}\au{u}^{6\gamma-2}|Au|^{4}.
  \end{eqnarray*}
  Using (\ref{stvro}) to estimate the last term in the previous 
  inequality, since $3p-2 \leq p$  finally we obtain
  \begin{eqnarray}
     \ep|\langle u'',\roep\rangle |(1+t)^{p} & \leq & 
     \frac{1}{4}\au{u}^{2\gamma}\auq{\roep}(1+t)^{p} + 
     c_{4}\ep^{2}\au{u}^{2(\gamma+1)} (1+t)^{p}+
      \nonumber \\
       &  &  +
       3\ep^{2}(1+t)^{5p}\au{u}^{6\gamma}|A^{3/2}u|^{2}.
      \label{spsi3-1}
  \end{eqnarray}
  From (\ref{spsi3-1}), (\ref{dv1}), (\ref{dv6}) thus we arrive at
  \begin{eqnarray}
      |\psi_{3}| &\leq  &
      \frac{1}{4}\int_{0}^{t}\au{u(s)}^{2\gamma}\auq{\roep(s)}(1+s)^{p}ds + 
      c_{5}\ep^{2} \nonumber \\
      & & \leq \frac{1}{4}\int_{0}^{t}
      \phi_{\rho}(s)(1+s)^{p}ds + 
      c_{5}\ep^{2}.
      \label{spsi3ro}
   \end{eqnarray}
  Let us now consider $\psi_{4}$ and prove that
  \begin{equation}
      p(1-p)\int_{0}^{t}(1+s)^{p-2}|\roep(s)|^{2}ds \leq c_{6}, \quad 
      \quad \forall t \geq 0.
      \label{spsi4ro-1}
  \end{equation}
  If $p=1$ thesis is obvious. If $p<1$ it is enough to prove that 
  $|\roep|^{2}$ is bounded independently from $\ep$ and $t$. But this 
  is a straightaway consequence of (\ref{D0}) and (\ref{dv1}). 
  
  Using (\ref{spsi4ro-1}) we then obtain
  \begin{equation}
      \psi_{4} \leq c_{6}\ep.
      \label{spsi4ro}
  \end{equation}
Now we go back to (\ref{stDro1-1}), and using (\ref{spsi1ro}), 
 (\ref{spsi3ro}), (\ref{spsi4ro}), since $p\ep \leq 1/4$ we achieve
 \begin{equation}
     \displaystyle \frac{1}{4}\int_{0}^{t}\phi_{\rho}(s)(1+s)^{p}ds   + 
       \frac{1}{8}|\roep(t)|^{2}\leq c_{7}\ep^{2}+c_{6}\ep+\ep 
       \displaystyle  \int_{0}^{t}|\roep'(s)|^{2}(1+s)^{p}ds.
     \label{stDro}
 \end{equation}
 Let us now remark that thanks to (\ref{D4}) and (\ref{dv5}) we have
 \begin{equation}
     \int_{0}^{t}|\roep'(s)|^{2}(1+s)^{p}ds \leq c_{8}.
     \label{SIro}
 \end{equation}
  Plugging  (\ref{SIro}) in (\ref{stDro}) we gain
  $$|\roep|^{2} \leq c_{9}\ep.$$
  At this point we can improve estimates (\ref{spsi4ro-1}) and 
  (\ref{spsi4ro}) as below
 $$  p(1-p)\int_{0}^{t}(1+s)^{p-2}|\roep(s)|^{2}ds \leq c_{10}\ep, 
       \hspace{1em} \psi_{4} \leq c_{10}\ep^{2}.
      $$
  Using this last estimate in (\ref{stDro}) instead of 
  (\ref{spsi4ro}) we finally get
  $$\displaystyle \frac{1}{4}\int_{0}^{t}\phi_{\rho}(s)(1+s)^{p}ds   + 
       \frac{1}{8}|\roep(t)|^{2}\leq c_{11}\ep^{2}+\ep 
       \displaystyle  \int_{0}^{t}|\roep'(s)|^{2}(1+s)^{p}ds.$$
  that is (\ref{stDro1}).
  \subparagraph{\textmd{\textsl{Proof of (\ref{stFro1})}}}
  From (\ref{dv4}) we have
  $$2\gamma \au{u}^{2(\gamma-1)}|Au|^{2}(1+t)^{p}\leq 
  \frac{2\gamma\gamma_{4}}{1+t}\leq 
  \frac{2\gamma\gamma_{4}}{(1+t)^{p}},$$
  hence from (\ref{hep1}), (\ref{DFro}), (\ref{sdiff}) and 
  (\ref{sv''}) we obtain
  \begin{eqnarray}
      F_{\rho}' & \leq & 
      -\frac{7}{4}\frac{|\rep'|^{2}}{\au{u}^{2\gamma}} 
      \frac{1}{(1+t)^{p}}+ 
      2\au{\roep}|A^{1/2}\tetep|+ \frac{1}{2}\frac{|\rep'|^{2}}{\au{u}^{2\gamma}} 
      \frac{1}{(1+t)^{p}}
      +
      \nonumber \\
       &  & +4\frac{(\au{u}^{2\gamma} - 
      \au{\uep}^{2\gamma})^{2}}{\au{u}^{2\gamma}}|A\uep|^{2}(1+t)^{p}+ 4\ep^{2}
      \frac{|u''|^{2}}{\au{u}^{2\gamma}} 
     (1+t)^{p}
      \nonumber  \\
       & \leq &  -\frac{5}{4}\frac{|\rep'|^{2}}{\au{u}^{2\gamma}} 
      \frac{1}{(1+t)^{p}}+ 
      2\au{\roep}|A^{1/2}\tetep|+
      \nonumber  \\
       &  & +24\gamma^{2} 
       \frac{|A\uep|^{2}(\au{u}^{2(\gamma-1)}+
     \au{\uep}^{2(\gamma-1)})}{\au{u}^{2\gamma}} \phi_{\rho}(1+t)^{p}+
      \nonumber  \\
       &  & + 4c_{1}\ep^{2}(1+t)^{3p-2}\au{u}^{2\gamma}|Au|^{2} 
    + 12\ep^{2}(1+t)^{5p}\au{u}^{6\gamma}|A^{2}u|^{2}.
      \label{sDFro-1} 
  \end{eqnarray}
  Let us now observe that, thanks to (\ref{D2}), (\ref{dv2}) and 
  (\ref{D1}), we have
 $$
      \frac{1}{\au{u}^{2\gamma}} |A\uep|^{2}\au{\uep}^{2(\gamma-1)}+
         \frac{1}{\au{u}^{2}}|A\uep|^{2}\au{\uep}^{2(\gamma-1)} 
	 \frac{1}{\au{\uep}^{2(\gamma-1)}}\leq c_{12}. 
      $$
  Replacing this inequality in (\ref{sDFro-1}) and integrating we get
  \begin{eqnarray*}
     &\displaystyle  F_{\rho}(t)+ \frac{5}{4}\int_{0}^{t}
     \frac{|\rep'(s)|^{2}}{\au{u(s)}^{2\gamma}} 
      \frac{1}{(1+s)^{p}} ds \leq  2\sup_{0\leq s\leq 
      t}\au{\roep(s)}\int_{0}^{t}|A^{1/2}\tetep(s)|\,ds+ & \\
       & \displaystyle   +c_{13}\int_{0}^{t}(1+s)^{p} \phi_{\rho}(s) ds+
       4c_{1}\ep^{2}\int_{0}^{t}(1+s)^{3p-2}\au{u(s)}^{2\gamma}|Au(s)|^{2}ds +& \\
       & \displaystyle  
       +
       12\ep^{2}\int_{0}^{t}(1+s)^{5p}\au{u(s)}^{6\gamma}|A^{2}u(s)|^{2}ds.& 
  \end{eqnarray*}
  We can now use Lemma \ref{thetalemma} with $\delta = 0$ and $j = 
  1$, (\ref{stDro1}), (\ref{dv0}) with $k=1$, (\ref{dv7}) and (\ref{dv1}), 
  thus we obtain
  \begin{eqnarray*}
     &\displaystyle  F_{\rho}(t)+ \frac{5}{4}\int_{0}^{t}
     \frac{|\rep'(s)|^{2}}{\au{u(s)}^{2\gamma}} 
      \frac{1}{(1+s)^{p}} ds    \leq  c_{14}\ep\sup_{0\leq s\leq 
      t}\au{\roep(s)} + c_{15}\ep^{2}+ & \\
       & \displaystyle+c_{16}\ep \int_{0}^{t}|\roep'(s)|^{2}(1+s)^{p}ds+ 
       4c_{1}\ep^{2}\frac{|Au_{0}|^{2}}
       {\auq{u_{0}}}\int_{0}^{t}(1+s)^{p}\au{u(s)}^{2(\gamma+1)}ds &   \\
       &\displaystyle \leq \frac{1}{2}\sup_{0\leq s\leq 
       t}\au{\roep(s)}^{2}+c_{17}\ep^{2}+
      c_{16}\ep \int_{0}^{t}|\roep'(s)|^{2}(1+s)^{p}ds.& 
  \end{eqnarray*}
  Let now $T>0$ and let us take the essup on $0\leq t\leq T$, then 
  we get
  $$\frac{1}{2}\sup_{0\leq t\leq T} F_{\rho}(t)+ 
  \frac{5}{4}\int_{0}^{T}
     \frac{|\rep'(s)|^{2}}{\au{u(s)}^{2\gamma}} 
      \frac{1}{(1+s)^{p}} ds \leq  c_{17}\ep^{2}+
      c_{16}\ep \int_{0}^{T}|\roep'(s)|^{2}(1+s)^{p}ds.$$
      Since $T$ is arbitrary we have proved (\ref{stFro1}).
     \subparagraph{\textmd{\textsl{Proof of (\ref{stEro1})}}}
      Let us estimate separately the terms $S_{1},\ldots, S_{5}$ in 
      (\ref{DEro}).
      
      Since $\ep\leq 1/4$ then
  \begin{equation}
      S_{1}\leq - \frac{3}{2}|\rep'|^{2}(1+t)^{p}.
      \label{stS1ro}
  \end{equation}
  Moreover
  \begin{equation}
      S_{2}\leq 2(1+t)^{p}{\au{u}^{2\gamma}}\auq{\roep} \leq 
      2(1+t)^{p}\phi_{\rho}.
      \label{stS2ro}
  \end{equation}
  Since $\au{u}^{2\gamma}$ is bounded then
  \begin{equation}
      S_{3}\leq c_{18}(1+t)^{2p}\au{\roep}\au{\tetep'}.
      \label{stS3ro}
  \end{equation}
  From (\ref{sv''}) we deduce
  \begin{eqnarray}
      S_{4}&\leq& 
      \frac{1}{4}|\rep'|^{2}(1+t)^{p}+4\ep^{2}|u''|^{2}(1+t)^{3p}\nonumber
      \\
      &\leq & \frac{1}{4}|\rep'|^{2}(1+t)^{p} 
      +4c_{1}\ep^{2}(1+t)^{3p}\au{u}^{4\gamma}|Au|^{2}+
      \nonumber \\
      & &+
      12\ep^{2}(1+t)^{7p}\au{u}^{8\gamma}|A^{2}u|^{2}.
      \label{stS4ro}
  \end{eqnarray}
  Let us now estimate $S_{5}$. From (\ref{sdiff}) and (\ref{D2}) we get
  \begin{eqnarray}
     & S_{5} \displaystyle \leq \frac{1}{4}|\rep'|^{2}(1+t)^{p} + 
      4(1+t)^{3p}|A\uep|^{2}(\au{u}^{2\gamma} - 
      \au{\uep}^{2\gamma})^{2}  &
      \nonumber  \\
       & \displaystyle \leq \frac{1}{4}|\rep'|^{2}(1+t)^{p} + 
       c_{19}(1+t)^{2p}|A\uep|^{2}\au{\uep}^{2(\gamma-1)}
       \left(1+\frac{\au{u}^{2(\gamma-1)}}{\au{\uep}^{2(\gamma-1)}}\right)
       \phi_{\rho}(1+t)^{p} &
      \nonumber  \\
       & \displaystyle \leq  \frac{1}{4}|\rep'|^{2}(1+t)^{p} 
       + c_{20}\frac{1}{(1+t)^{1-p}}
       \left(1+\frac{\au{u}^{2(\gamma-1)}}{\au{\uep}^{2(\gamma-1)}}\right) 
       \phi_{\rho}(1+t)^{p}.&
      \label{stS5ro-1}
  \end{eqnarray}
  Now we want to  prove that 
  \begin{equation}
      \chi:= \frac{1}{(1+t)^{1-p}}\frac{\au{u}^{2(\gamma-1)}}{\au{\uep}^{2(\gamma-1)}} \leq c_{21}.
      \label{stS5ro-2}
  \end{equation}
  When $A$ is coercive this is a consequence of (\ref{D1}) and 
  (\ref{dv3}). On the other hand if $p\leq 
  (\gamma^{2}+1)/(\gamma^{2}+2\gamma -1)$ then 
  $$ \alpha: = \frac{\gamma^{2}
  +1-p(\gamma^{2}+2\gamma-1)}{\gamma(\gamma+1)} \geq 0,
  $$
  hence
  from (\ref{D1}) and 
  (\ref{dv2}) we deduce:
  $$ \chi
  \leq \frac{c_{22}}{(1+t)^{1-p}} \frac{(1+t)^{(p+1)
  (\gamma-1)/\gamma}}{(1+t)^{(p+1)(\gamma-1)/(\gamma+1)}}= 
  \frac{c_{22}}{(1+t)^{\alpha}}\leq c_{22}.$$
  At this point from (\ref{stS5ro-1}) and (\ref{stS5ro-2}) it follows 
  that:
  \begin{equation}
      S_{5}\leq  \frac{1}{4}|\rep'|^{2}(1+t)^{p} 
       +  
       c_{23}\phi_{\rho}(1+t)^{p}.
      \label{stS5ro}
  \end{equation}
  If we put (\ref{stS1ro}), (\ref{stS2ro}), (\ref{stS3ro}), 
  (\ref{stS4ro}), (\ref{stS5ro}) in (\ref{DEro}), since $S_{6} \leq 0$ we 
  get:
  \begin{eqnarray*}
      E_{\rho}' + |\rep'|^{2}(1+t)^{p}  & \leq & 
      c_{24}\phi_{\rho}(1+t)^{p}+ c_{18}(1+t)^{2p}\au{\roep}\au{\tetep'} +  \\
       &  & +
       4c_{1}\ep^{2}(1+t)^{3p}\au{u}^{4\gamma}|Au|^{2} + 12\ep^{2}(1+t)^{7p}\au{u}^{8\gamma}|A^{2}u|^{2}.
  \end{eqnarray*}
  Integrating in $[0,t]$ and using (\ref{stDro1}), (\ref{dv5}), (\ref{dv7}) we 
  then achieve:
  \begin{eqnarray}
    & \displaystyle  E_{\rho}(t) + \int_{0}^{t}|\rep'(s)|^{2}(1+s)^{p} ds  \leq  
      c_{25}\ep^{2}+c_{26}\ep\int_{0}^{t}|\roep'(s)|^{2}(1+s)^{p}ds +&
      \nonumber  \\
       &  \displaystyle  +c_{18}\sup_{0\leq s \leq 
       t}(1+s)^{p}\au{\roep(s)}\au{u(s)}^{\gamma}
       \int_{0}^{t}\frac{\au{\tetep'(s)}}{\au{u(s)}^{\gamma}}(1+s)^{p}ds.&
      \label{stEro1-1}
  \end{eqnarray}
  Thanks to (\ref{dv2}) and  Lemma \ref{thetalemma}, with $j = 1$ and 
  $\delta = (3p+1)/2$ (note that thanks to (\ref{hep1}) all hypotheses of Lemma 
  \ref{thetalemma} are verified) we have
$$
     \int_{0}^{t}\frac{\au{\tetep'(s)}}{\au{u(s)}^{\gamma}}(1+s)^{p}ds 
     \leq c_{27}\int_{0}^{t}\au{\tetep'(s)}(1+s)^{(3p+1)/2}ds 
     \leq c_{28}\ep.
     $$
     Finally plugging this inequality in (\ref{stEro1-1}) and using 
   (\ref{SIro}) we gain
   \begin{eqnarray*}
        E_{\rho}(t) + \int_{0}^{t}|\rep'(s)|^{2}(1+s)^{p} ds  & \leq 
	& c_{25}\ep^{2}+c_{29}\ep +  c_{30}\ep\sup_{0\leq s \leq 
       t}(1+s)^{p}\au{\roep(s)}\au{u(s)}^{\gamma} \\
        &  & \leq c_{31}\ep + \frac{1}{2}\sup_{0\leq s \leq 
       t}E_{\rho}(s).
   \end{eqnarray*}
   Let now $T > 0$, if we take the essup for $0\leq t \leq T$ we get
   $$\frac{1}{2}\sup_{0\leq t \leq 
       T}E_{\rho}(t) +  \int_{0}^{T}|\rep'(s)|^{2}(1+s)^{p} ds \leq 
       c_{31}\ep.$$
        Since $T$ is arbitrary we have proved (\ref{stEro1}).
	\qed
\subsubsection{Conclusion}
From (\ref{SIro}), (\ref{stFro1}), (\ref{stDro1}) we straight obtain 
(\ref{H1}).

Let us now prove (\ref{H1c}). From (\ref{stEro1}) and Lemma 
\ref{thetalemma} with $j=0$ and $\delta = p$ we get
\begin{eqnarray*}
    \int_{0}^{+\infty}|\roep'(s)|^{2}(1+s)^{p} ds & \leq  &
    2\int_{0}^{+\infty}|\rep'(s)|^{2}(1+s)^{p}ds + 
    2\int_{0}^{+\infty}|\tetep'(s)|^{2}(1+s)^{p}ds
    \nonumber  \\
     & \leq & 2\gamma_{9}\ep + 2\ep C_{p} \sup_{s\geq 0}|\tetep'(s)|
     \leq C \ep,
    \label{SIro1}
\end{eqnarray*}
where $C$ depends only on the data, since $|\tetep'|\leq |w_{0}|$.
Finally replacing this estimate in (\ref{stDro1}) and (\ref{stFro1}) we obtain 
inequality (\ref{H1c}).

\label{NumeroPagine}

\end{document}